\documentclass[11pt]{article}
\usepackage{amsthm} 
\usepackage{amssymb, latexsym}
\usepackage{amsmath, amsfonts}

\usepackage{wrapfig}

\usepackage{verbatim}
\usepackage{epsfig}
\usepackage{epstopdf}
\usepackage{graphics,graphicx}
\usepackage{siunitx}
\usepackage[ruled,vlined]{algorithm2e}
\usepackage{array}
\usepackage{tabu}
\usepackage{color}
\usepackage{caption}
\usepackage{subcaption}
\usepackage{placeins}
\usepackage{float}
\usepackage{cleveref}

\usepackage{tikz}
\usepackage{wrapfig}
\usetikzlibrary{positioning} 
\usepackage{booktabs} 
\usepackage{multirow}

\usepackage{tikz-3dplot}
\usetikzlibrary{shapes}
\usepackage{pgfplots}
\usepackage{tikz}
\usetikzlibrary{intersections, pgfplots.fillbetween}
\usetikzlibrary{shapes}

\usepackage[left=2cm, right=2cm, top=2cm]{geometry}
\usepackage{pgfplots}
\usepackage{filecontents}

\newtheorem{theorem}{Theorem}[section]

\newtheorem{remark}[theorem]{Remark}

\newtheorem{assumption}[theorem]{Assumption}

\newcommand {\real} {\,\hbox{{\rm R \kern-1.2em I \kern.3em }}}

\def \bU{{\mathbf U}}

\def \bv{{\mathbf v}}
\def \b0{{\mathbf 0}}
\def \bu{{\mathbf u}}

\def \bbf{{\mathbf f}}
\def \bH{{\mathbf H}}

\def \bn{{\mathbf n}}

\def \bv{{\mathbf v}}
\def \bf{{\mathbf f}}

\def \bmu{\mbox{\boldmath $\mu$}}

\def \bSigma{\mbox{\boldmath $\Sigma$}}

\def \nn{\nonumber}

\def \Omegn1{\Omega^f_{t_{n+1}}}

\def \beq{\begin{equation}}
\def \eeq{\end{equation}}
\def \beqn{\begin{eqnarray}}
\def \eeqn{\end{eqnarray}}

{\color{red}
\title{\textbf{A novel approach to study the wellposedness of the 3D fluid-2D plate interaction PDE System}}}
\author{
George Avalos
\thanks{Department of Mathematics, University of Nebraska–Lincoln, Lincoln, NE 68588, USA ({\tt gavalos2@unl.edu}). 
}
\and 
Pelin G. Geredeli
\thanks{School of Mathematical and Statistical Sciences, Clemson
        University, Clemson, SC  29634, USA ({\tt pgerede@clemson.edu}). 
}
\and 
Hemanta Kunwar
\thanks{Department of Mathematics, Embry-Riddle Aeronautical University, Daytona Beach, FL 32114, USA ({\tt kunwarh@erau.edu}). 
}
\and 
Hyesuk Lee
  \thanks{School of Mathematical and Statistical Sciences, Clemson
        University, Clemson, SC  29634, USA ({\tt hklee@clemson.edu}). 
}
}
\date{}
\begin{document}

\maketitle 

\setcounter{equation}{0}
\setcounter{figure}{0}
\begin{abstract}
We consider a certain fluid-structure interaction (FSI) system with a view of obtaining an alternative methodology for establishing its strongly continuous semigroup wellposedness. (Semigroup generation for this FSI was originally considered in \cite{A-Clark}.) The FSI model under consideration describes the vibrations of an incompressible fluid within a 3D cavity as it interacts with the elastic membrane on the ``free" upper boundary of the cavity. Such coupled PDE systems appear in variety of natural settings such as biomedicine, aeroelasticity, and fluid dynamics.

Our proof of $C
_0$-semigroup wellposedness is based on a proper application of Lumer Phillips Theorem. In this regard, our main challenge is to show the maximality of the corresponding semigroup generator. To this end, we develop a ``nonstandard" inf-sup approach which avoids the use of technical nonlocal maps in the associated bilinear forms--unlike the earlier paper \cite{A-Clark}--and  allows for the solution of the fluid and plate solution variables simultanously. Our new inf-sup strategy will lead to a more efficient mixed finite element method (FEM) for approximating solutions
to the FSI problem, inasmuch our novel variational formulation avoids bilinear forms which are free from the computationally-intensive nonlocal solution operators invoked in \cite{A-Clark}.  We also perform numerical tests based on this formulation using a benchmark problem
 and present numerical results to demonstrate the effectiveness of our approach.
\end{abstract}

\begin{paragraph}{Keywords:}
{Stokes Fluid-Plate Interaction PDE System, Biharmonic Equation, Kirchoff Plate, wellposedness, domain decomposition.}
\end{paragraph}

\begin{paragraph}{2020  Mathematics Subject Classification:} 
{
Primary: 35Q35, Secondary: 65N30, 35A01, 74K20, 76M10
}
\end{paragraph}
\section{Introduction}

Coupled PDE systems which mathematically describe the interaction of a fluid flow with a
deformable elastic structure arise in various aeroelastic and biomedical applications; e.g., the fluttering of an airplane/airfoil wing, the movement of wind
turbines and bridges, blood transportation processes within arterial walls, the control of ocular pressure in the eye, and sloshing \cite{AD, FSIforBIO, Barbu, buk, Chambolle, hsu, Bociu, dow, eye}. As an exemplifying model, the fluid-structure interaction (FSI) system (\ref{1})-\eqref{4} below describes the vibrations of a fluid within a 3D cavity as it interacts with the elastic membrane on the ``free" upper boundary of the cavity. Since these systems are governed by multi-physics PDE\
models of generally high complexity, the theoretical and numerical analysis of solutions to such systems tend to be on the  challenging side.

Moreover, although there have been extensive studies on numerical methods for standard FSI problems with moving or {{non-moving}} interfaces, where the structure is modeled by linear {\cite{lin1,lin2} or nonlinear elasticity \cite{Barbu, nl2, nl3}}--such problems are broadly classified as monolithic or partitioned methods-- the body of work with respect to the quantitative properties and numerical analysis of 3D Stokes/Navier Stokes flow and 2D plate interaction FSI systems is not extensive \cite{Barbu,lin2,lin1,nl3, A-Clark, P8}. 

From the theoretical point of view, the wellposedness of the particular FSI model in (\ref{1})-\eqref{4} was initially shown in \cite{CR, A-Clark}. While \cite{CR} employed a Galerkin approximation for the Stokes-nonlinear plate interaction FSI model in which nonlinearity appears strictly in the plate component, \cite{A-Clark} invoked a strongly continuous semigroup approach for the fully linear FSI system. In \cite{A-Clark}, the particular wellposedness argument ultimately provides a mixed finite element method (FEM) formulation in order to approximate the solution of the FSI system (\ref{1})-\eqref{4}.  In particular, the inf-sup formulation in \cite{A-Clark} is predicated on recovering only the elastic solution component of the FSI system. The bilinear forms which are created for this (structural) mixed variational formulation of \cite{A-Clark} involve nonlocal (fluid) solution maps which essentially describe the coupling of the two dimensional  plate with the three dimensional Stokes dynamics, and are unavoidably complicated in form. Upon solving the plate inf-sup problem, the fluid solution component of the FSI is subsequently recovered in \cite{A-Clark} via said nonlocal maps. A companion mixed finite element method (FEM) naturally arises, which is based on the (plate) inf-sup approach in \cite{A-Clark}. However, its implementation becomes quite computationally intensive, even for the static case. 

By contrast, in this manuscript, our main goal is to simultaneously solve the fluid and plate solution variables via an inf-sup approach that appears ``natural'' to the given FSI. That is, there are no nonlocal maps in the associated bilinear forms; they naturally appear via standard invocations of Green's formulae. Our reasoning is that this wellposedness argument would in principle give rise to a more efficient mixed FEM scheme so as to approximate solutions of the FSI model, inasmuch as this FEM would deal with a variational formulation much simpler in structure than that created in \cite{A-Clark}. In order to generate this nonstandard variational formulation, the key point is to introduce a Lagrange multiplier $\mathfrak{g} \in   H^{-1/2}(\Omega_p)$ which is characterized by the boundary trace of the pressure, $p|_{\Omega_p}$ on $\Omega_p$. 

In this manuscript, we first establish, by our new approach, the strongly continuous semigroup wellposedness of the FSI model (\ref{1})-\eqref{4}. In particular, we associate solutions of the FSI dynamics with a $C_0$-semigroup of contractions. In this connection, our major challenge becomes establishing the maximality of the candidate semigroup generator via a \textit{nonstandard} mixed variational formulation which is wholly distinct from that of \cite{A-Clark}.

 Secondly, we study finite element approximation of the FSI model based on the mixed variational formulation, using the same numerical example reported in \cite{A-Clark}.  
 The example is a manufactured problem with a known exact solution, which allows us to test the convergence behavior of the finite element solutions. In this part of the work, we employ $P2-$Morley elements \cite{M4} for the plate subproblem and Taylor-Hood elements for the fluid subproblem. Numerical methods for fluid-structure interaction problems are typically classified as either monolithic or partitioned approaches. In the monolithic framework, the fluid and structural equations are combined into a single system and solved simultaneously, offering strong stability advantages in many scenarios. However, this strategy may be less suitable for fluid-plate interaction models, where the plate dynamics are governed by a fourth-order partial differential equation. In such cases, the differing regularity requirements between the fluid and plate subproblems pose significant difficulties. The main challenges associated with our analysis and the novelties in this manuscript are as follows:\\

 \newpage 

\textbf{1) Formulation of a nonstandard mixed variational problem.} In our work to show strongly continuos semigroup wellposedness for the FSI system (\ref{1})-\eqref{4}, the proof of maximality of the semigroup generator is the crux of the matter. For this, it is necesssary to consider a static (essentially Laplace-transformed) version of PDE system (\ref{1})-\eqref{4}. A multiplication of the governing static FSI equations by appropriate test functions will ultimately give rise to a certain saddle point problem whose solution relies on a Babuska-Brezzi (inf-sup) approach. At this point, we emphasize that because of the matching fluid and structure velocities, as well as the presence of the fluid pressure as an external source in the plate dynamics, the argument we follow here for the necessary Babuska-Brezzi (inf-sup) formulation is very different than that given in \cite{A-Clark} or for other FSI. In particular, our mixed variational formulation is created with respect to \textit{all} solution components $\{u,w,p,p|_{\Omega_p}\}$. \\

\textbf{2) Establishment of the Inf-sup Condition.} Our mixed variational formulation mentioned above in (1) is given via the bilinear forms $a_{\lambda}(.,.)$ (that are composed of bilinear forms which are ``natural'' to fluid and elastic solution variables), and $b(.,.)$ which accounts in part for the underlying constraints coming from the matching velocities boundary conditions and the pressure term on $\Omega_p$-- it is not solely the usual divergence free constraint for incompressible flows (see \eqref{a-l} and \eqref{eq75}).  With  $\Sigma \equiv \mathbf{U}\times H_{0}^{2}(\Omega _{p})$, where $%
\mathbf{U}\subset \mathbf{H}^{1}(\Omega _{f})$, and $L_{0}^{2}(\Omega _{f})$ is
the zero average value subspace of $L^{2}(\Omega _{f})$,
our major challenge will be to show the inf-sup condition for the bilinear form $b(.,.).$ to establish the inequality for some $\beta> 0:$ $$sup_{(\bv,z)\in \Sigma}\frac{b((\bv,z), (l,h,r))}{||(\bv,z)||_{\Sigma}}\geq \beta (||l||_{\Omega_p}+||h||_{[H^{1/2}_{\overline{\Omega}}]'}+|r|),~~~\forall~(l,h,r)\in L_0^2(\mathcal{O})\times[H^{1/2}_{\overline{\Omega}}]'\times\mathbb{R},$$ where the introduction of the space $H^{1/2}_{\overline{\Omega}}=\{\eta \in H^{1/2}(\mathbb{R}^2): supp ~\eta \subset \overline{\Omega}\}$ is very crucial in our pending analysis, to extend certain $H^{1/2}(\Omega_p)$-boundary data by zero outside $\Omega_p$ and still maintain $H^{1/2}$-regularity. \\

\textbf{3) Domain Decomposition for Numerical Simulation.} The biharmonic operator which dictates the plate component of the FSI presents a numerical challenge  due to the requirement of the use
of class $C^1$-class  finite elements for accurate approximation, whereas the fluid velocity can be simulated using $C^0$-class elements. To address this in our numerical example below, we adopt a domain decomposition approach for the coupled PDE system and design a decoupling algorithm that allows the fluid and plate subproblems to be solved independently. 
This algorithm is an implicit type for better stability, but it requires iterations between subproblems. 
 The regularity challenges introduced by the biharmonic operator in the plate equation can be tackled through various strategies, including mixed finite element methods \cite{mix1,mix2,mix3}, high-order conforming elements \cite{High1,High2}, or non-conforming elements like the Morley element \cite{M1,M2,M3}.

\vspace{0.5cm}

\textbf{Plan of the paper.} The remainder of the paper is structured as follows. In \Cref{meq}, we present the 3D fluid–2D plate interaction model along with the notations used throughout the work. The main well-posedness result for the coupled system is established in \Cref{sec:wellposedness}. Next, we outline the numerical algorithm and showcase the corresponding results in \Cref{sec:NR}. Finally, concluding remarks are provided in \Cref{sec:Conclude}.

\section{Description of the PDE Model}\label{meq}
In this manuscript, the geometry of the coupled PDE dynamics will be designated as follows: the structure domain -- i.e., the portion of the geometry in which the structural PDE component evolves -- is denoted by $\Omega_p$, where $\Omega_p \subset \{x\in \mathbb{R}^3: x=(x_1,x_2,0)\}$. Moreover, $S$ will be a surface contained in  $ \left\{ x\in \mathbb{R}^3: x=(x_{1},x_{2},x_{3}):x_{3}\leq 0\right\}$. Therewith, the fluid domain $\Omega_f \subset \real^3$ is so configured that its boundary satisfies $\partial \Omega_f := \overline{\Omega_p} \cup \overline{S}$. Thus, the fluid domain $\Omega_f$ is occupied by the free fluid while the domain $\Omega_p$ constitutes the plate region.
\begin{figure}[H]
\centering
\begin{tikzpicture}[scale=0.9]
\draw[left color=black!10,right color=black!20,middle
color=black!50, ultra thick] (-2,0,0) to [out=0, in=180] (2,0,0)
to [out=270, in = 0] (0,-3,0) to [out=180, in =270] (-2,0,0);

\draw [fill=black!60, ultra thick] (-2,0,0) to [out=80,
in=205](-1.214,.607,0) to [out=25, in=180](0.2,.8,0) to [out=0,
in=155] (1.614,.507,0) to [out=335, in=100](2,0,0) to [out=270,
in=25] (1.214,-.507,0) to [out=205, in=0](-0.2,-.8,0) [out=180,
in=335] to (-1.614,-.607,0) to [out=155, in=260] (-2,0,0);

\draw [dashed, thin] (-1.7,-1.7,0) to [out=80, in=225](-.6,-1.3,0)
to [out=25, in=180](0.35,-1.1,0) to [out=0, in=155] (1.3,-1.4,0)
to [out=335, in=100](1.65,-1.7,0) to [out=270, in=25] (0.9,-2.0,0)
to [out=205, in=0](-0.2,-2.2,0) [out=180, in=335] to (-1.514,-2.0)
to [out=155, in=290] (-1.65,-1.7,0);

\node at (0.2,0.1,0) {{\LARGE$\Omega_p$}};

\node at (1.95,-1.5,0) {{\LARGE $S$}};

\node at (-0.3,-1.6,0) {{\LARGE $\Omega_f$}};
\end{tikzpicture}

\caption{Fluid-Structure Interaction Geometry }
\end{figure}
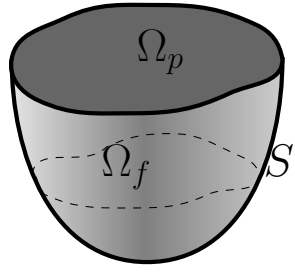
Posed on finite energy space $\mathbf{H}$, to be explicitly identified below, we consider the flow in $\Omega_f$ to be governed by the time-dependent Stokes equation and the plate structure is modeled by the ``Euler-Bernoulli” or ``Kirchhoff" equation. The fluid-plate system is then written as, for given terminal time $T > 0$:
\begin{equation}
\label{1}
\begin{cases}
\bu_t - \nu_f \Delta \bu + \nabla p = \mathbf{f}_f & \text{in} \ \Omega_f \times (0, T), \\
\nabla \cdot \bu = 0~~ \text{in}~~ \ \Omega_f \times (0, T), 
\end{cases}
\end{equation}
\begin{align}
\label{2}
\begin{cases}
& w_{tt} - \rho \Delta w_{tt} + \Delta^2 w =p \big|_{\Omega_p} +  f_p  \quad \text{in} \ \Omega_p \times (0, T), \\
& w = \frac{\partial w}{\partial \bn_p} = 0 \quad \text{on} \ \partial \Omega_p,
\end{cases}
\end{align}

\begin{equation}
\mathbf{u}=[u_{1},u_{2},u_{3}]=\left\{ 
\begin{array}{l}
\mathbf{0}~~\text{, \ on }S\times (0,T) \\ 
\lbrack 0,0,w_{t}]~~\text{, \ on }\Omega _{p}\times (0,T)%
\end{array}%
\right.   \label{3}
\end{equation}

 \begin{align}
 \label{4}
 [u(0),w(0),w_t(0)]=[u_0,w_0,w_{1}] \in \mathbf{H}.
\end{align}
Here, solution variables $\bu(\textbf{x},t)$ and $p(\textbf{x},t)$ denote the fluid velocity and the fluid pressure, respectively in $\Omega_f$. Moreover, $\nu_f$  denotes the constant fluid viscosity, and $\bbf_f(t) \in \mathbf{L}^{2}(\Omega _{f})$ and $f_p(t)  \in L^{2}(\Omega _{p})$ are given forces . Furthermore, $w(\textbf{x},t)$ is the displacement of the plate structure; $\bn_p$ denotes the outward unit normal 
vector to $\Omega_p$. The constant $\rho \geq 0$ is the rotational inertia parameter which is proportional to the square of the thickness of the plate. The case $\rho
=0$ indicates that there is no rotational plate dynamics in play (i.e., the thickness of the plate is neglected), whereas 
$\rho >0$ refers to the presence of rotational forces (that is, the full Kirchoff plate).   Also, $\mathbf{\nu }(x)$ is the unit outward normal vector with respect to $%
\Omega _{f}$.
  \\

\noindent Since the geometry described above is non-smooth, we will eventually require the following higher regularity results for FSI solutions. To have the validity of these regularity results, the geometry $\{\Omega_f,\Omega_p\}$ is assumed to satisfy the following \cite{avalos2014exponential}:
 \begin{assumption} \label{geometry}
 With respect to the above geometry, the pair $[\Omega_f, \Omega_p]$ satisfies one of the following assumptions
\begin{enumerate}
    \item $\Omega_f$ is a convex domain with wedge angles $\leq \frac{2\pi}{3}$. Moreover, $\Omega_p$ has smooth ($C^4$) boundary, and $S$ is a piecewise smooth surface.
    \item $\Omega_f$ is convex polyhedron having angles $\leq \frac{2\pi}{3}$, and so then $\Omega_p$ is a convex polygon with angles $\leq \frac{2\pi}{3}$.
\end{enumerate}
  \end{assumption}
  \vspace{0.2cm}
 
 \noindent As shown in \cite{avalos2014exponential}, under Assumption \ref{geometry}, solutions of the FSI system (\ref{1})-(\ref{4}) with initial data in the domain of the associated semigroup generator -- described explicitly in \cite{A-Clark} -- satisfy pointwise in time,
\[
\left\{ \mathbf{u},p\right\} \in \mathbf{H}^{2}(\Omega _{f})\times
H^{1}(\Omega _{f})\text{ and }\left[ w,w_{t}\right] \in \left[ H^{s_{\rho
}}(\Omega _{p})\cap H_{0}^{2}(\Omega _{p})\right] \times H_{0}^{2}(\Omega
_{p}),\text{where }s_{\rho }=\left\{ 
\begin{array}{c}
4\text{, if }\rho =0 \\ 
3\text{, if }\rho >0%
\end{array}%
\right. .
\]

 \noindent \textbf{Notation}
   \vspace{0.2cm}
   
\noindent  The space $H^{s}(D)$ will denote the Sobolev
space of order $s$, defined on a domain $D$, and $H_{0}^{s}(D)$ denotes the
closure of $C_{0}^{\infty }(D)$ in the $H^{s}(D)$ norm which we denote by $%
\Vert \cdot \Vert _{H^{s}(D)}$ or $\Vert \cdot \Vert _{s,D}$. Norms $||\cdot ||_{D}$ are taken to be $%
L^{2}(D)$ for the domain $D$. Inner products in $L^{2}(D)$ are written $%
(\cdot ,\cdot )_{D}$, and the dual product between $H^s$ and its dual space is denoted by %
$\langle \cdot ,\cdot \rangle_{D} $.\\

\noindent In formulating the appropriate finite energy space $\mathbf{H}$ for the FSI system (\ref{1})-(\ref{4}), we must account for an underlying compatibility condition satisfied by the
structural PDE component $[w,w_t]$. Indeed, invoking the
matching velocities boundary condition in (\ref{3}), the divergence free
constraint on $\mathbf{u}$, and the fact that normal vector $\mathbf{\nu }%
=[0,0,1]$ on $\Omega _{p}$, we have by Green's Identity,%
\begin{equation}
\begin{array}{l}
\frac{d}{dt}\int_{\Omega _{p}}w(t)d\Omega _{p}=\int_{\Omega
_{p}}w_{t}(t)d\Omega _{p}=\int_{\Omega _{p}}u_{3}(t)d\Omega
_{p}=\int_{\Omega _{f}}\nabla \cdot \mathbf{u}d\Omega _{f}=0, \\ 
\\ 
\text{and so, pointwise in time} \\ 
\\ 
\int_{\Omega _{p}}w(t)d\Omega _{p}=\text{ \emph{constant}.}%
\end{array}
\label{compat}
\end{equation}

\noindent With this constraint in mind, we proceed to construct $\mathbf{H}$ and other needed energy spaces:
 \begin{align}
&\bU :=  \{\bv = (v_1,v_2,v_3) \in \bH^{1}(\Omega_f):\; v_1=v_2=0 \mbox{ on } \Omega_p, \bv = \textbf{0} \mbox{ on } S \}. \label{U} 
\end{align}
We also define the mean zero space for the pressure,
\begin{equation}
Q\equiv \{q\in L^{2}(\Omega _{f}):\int_{\Omega _{f}}q\ d\Omega _{f}=0\}.
\label{Q}
\end{equation}
Moreover, in what follows, we will be under the necessity of extending boundary data on $\Omega_p$ by zero and maintaining full Sobolev regularity on $\partial \Omega _{f}$, which is not possible for a general $H^{1/2}(\Omega_p)$-function (see e.g., Theorem 3.33, p. 95 of \cite{McLean}). For this reason, we need the following space, as well as its topological dual: 
\begin{equation}
H_{\overline{\Omega }_{p}}^{\frac{1}{2}}=\left\{ z\in H^{\frac{1}{2}}(%
\mathbb{R}^{2}):\text{supp}(z)\subset \overline{\Omega }_{p}\right\} .
\label{contain}
\end{equation}
Bearing in mind the compatibility condition in (\ref{compat}) for the structural PDE component of FSI, the energy space $\mathbf{H}$, associated with (\ref{1})-(\ref{4}) is defined as follows: 
Let 
\begin{equation}
\begin{array}{c}
\mathcal{H}_{\text{fluid}}\equiv \left\{ \mathbf{f}\in \mathbf{L}^{2}(\Omega
_{f}):\text{ }\nabla \cdot \mathbf{f}=0\text{, }\left[ \mathbf{f}\cdot 
\mathbf{\nu }\right] _{S}=0\right\} ; \\ 
\\ 
\widehat{L}^{2}(\Omega _{p})=\left\{ w\in L^{2}(\Omega _{p}):\int_{\Omega
_{p}}w~d\Omega _{p}=0\right\} ;%
\end{array}
\label{energy_1}
\end{equation}

and
\begin{equation}
W_{\rho }=\left\{ 
\begin{array}{l}
\widehat{L}^{2}(\Omega _{p})\text{, if }\rho =0 \\ 
H_{0}^{1}(\Omega _{p})  \cap \widehat{L}^{2}(\Omega _{p})\text{, if }\rho >0%
\end{array}%
\right. .  \label{energy_2}
\end{equation}

Therewith, we then set
\begin{equation}
    \bH = \{[\mathbf{f}, \omega_1, \omega_2] \in\mathcal{H}_{\mbox{fluid}} \times  [H_0^{2}(\Omega_p) \cap \widehat{L}^{2}(\Omega_p)] \times W_{\rho}  \mbox{ with } \mathbf{f} \cdot \nu|_{\Omega_p} = \omega_2 \}.
\end{equation}
%
 \section{Preliminaries}

\noindent Our methodology to show the existence-uniqueness of the solution to the PDE system \eqref{1}-\eqref{4} is based on the application of Lumer Philips Theorem. Hence, the elimination of the pressure variable \cite{A-Clark} is very crucial in order to formulate the PDE system \eqref{1}-\eqref{4} as an abstract ODE in Hilbert Space $\mathbf{H}.$ For the readers' convenience, we provide here the particulars in \cite{A-Clark} concerning the associated FSI $C_0$- semigroup generator.
\bigskip 

\noindent We start by defining the following realization of the translated Laplacian:
Given $\rho \geq 0$, we define the positive definite, self-adjoint operator $%
P_{\rho }:L^{2}(\Omega _{\rho })\rightarrow $ $L^{2}(\Omega _{\rho })$ via %
\begin{equation}
P_{\rho }\varpi =\varpi -\rho \Delta \varpi \text{, \ for }\varpi \in
D(P_{\rho })=\left\{ 
\begin{array}{l}
L^{2}(\Omega _{\rho })\text{, if }\rho =0 \\ 
H^{2}(\Omega _{p})\cap H_{0}^{1}(\Omega _{p})\text{, if }\rho >0%
\end{array}%
\right. .  \label{P}
\end{equation}
Therewith, an application of the divergence operator to the Stokes PDE component in \eqref{1}-\eqref{4}, and the fact that $\bu$ is solenoidal gives that the associated pressure variable $p(t)$ is harmonic pointwise in time, and in fact evidently solves the BVP,
\begin{align}
& \Delta p=0 \text { in } \Omega_f  \label{p-12}\\
& \frac{\partial p}{\partial \nu}+P_{\rho}^{-1} p=P_{\rho}^{-1} \Delta^{2} w+\left.\Delta u^{3}\right|_{\Omega} \quad \text { on } \Omega_p  \label{p-13}\\
& \frac{\partial p}{\partial \nu}=\left.\Delta \bu \cdot \nu\right|_{S} \text { on } S. \label{p-14}
\end{align}
In generating the boundary condition in \eqref{p-14}, we formally multiplied both sides of the fluid equation in \eqref{1} with an appropriate extension (into $\Omega_f$) of unit normal vector $\nu$, and then subsequently take the resulting trace on $S$. (We implicitly used the fact that $\bu=0$ on $S$.) Also, in establishing the boundary condition in \eqref{p-13}, we used the fact that $\nu=[0,0,1]$ on $\Omega_p$, the plate equation in \eqref{2} and the boundary condition \eqref{3}, to have
$$
\begin{aligned}
P_{\rho}^{-1} \Delta^{2} w & =-w_{t t}+\left.P_{\rho}^{-1} p\right|_{\Omega_p} \\
& =-\frac{d}{d t}\left[0,0, w_{t}\right] \cdot \nu+\left.P_{\rho}^{-1} p\right|_{\Omega_p} \\
& =-\left[\bu_{t} \cdot \nu\right]_{\Omega_p}+\left.P_{\rho}^{-1} p\right|_{\Omega_p} \\
& =-[\Delta \bu \cdot \nu]_{\Omega_p}+\left.\frac{\partial p}{\partial \nu}\right|_{\Omega}+\left.P_{\rho}^{-1} p\right|_{\Omega_p},
\end{aligned}
$$
which gives \eqref{p-13}. We note that the BVP (\ref{p-12})-(\ref{p-14}) can be solved through the agency of the following ``Robin" maps $R_{\rho}$ and $\tilde{R}_{\rho}$ : We define
\begin{align*}
& R_{\rho} g=f \Leftrightarrow\left\{\Delta f=0 \text { in } \Omega_f ;~~ \frac{\partial f}{\partial \nu}+P_{\rho}^{-1} f=g \text { on } \Omega_p ;~~ \frac{\partial f}{\partial \nu}=0 \text { on } S\right\} ; \\
& \tilde{R}_{\rho} g=f \Leftrightarrow\left\{\Delta f=0 \text { in } \Omega_f ;~~ \frac{\partial f}{\partial \nu}+P_{\rho}^{-1} f=0 \text { on } \Omega_p ;~~ \frac{\partial f}{\partial \nu}=g \text { on } S\right\} . 
\end{align*}
By the Lax-Milgram Theorem we have,
\begin{equation*}
R_{\rho} \in \mathcal{L}\left(H^{-\frac{1}{2}}(\Omega_p), H^{1}(\Omega_f )\right) ; \quad \tilde{R}_{\rho} \in \mathcal{L}\left(H^{-\frac{1}{2}}(S), H^{1}(\Omega_f )\right) .
\end{equation*}
Therefore, the pressure variable $p(t)$, as necessarily the solution of  (\ref{p-12})-(\ref{p-14}), can be written pointwise in time as
\begin{equation}
\label{pressure}
p(t)=G_{\rho, 1}(w(t))+G_{\rho, 2}(\bu(t)),
\end{equation}
where
\begin{align}
G_{\rho, 1}(w) & =R_{\rho}\left(P_{\rho}^{-1} \Delta^{2} w\right)   \label{G1} \\
G_{\rho, 2}(\bu) & =R_{\rho}\left(\left.\Delta u^{3}\right|_{\Omega_p}\right)+\tilde{R}_{\rho}\left(\left.\Delta \bu \cdot \nu\right|_{S}\right) . \label{G2}  
\end{align}
In short, the formal vector calculations above which give rise to the BVP (\ref{p-12})-(\ref{p-14}) allow for an elimination of the pressure variable $p(t)$ in \eqref{1}-\eqref{4}. In consequence, and as shown in \cite{A-Clark}, the PDE system \eqref{1}-\eqref{4} is written as the abstract Cauchy problem,
\begin{equation}
\label{ODE}
   \begin{cases} \frac{d}{dt} \begin{bmatrix}
w(t)\\
w_t(t)\\
\bu(t)
\end{bmatrix} = \mathbf{A}_\rho \begin{bmatrix}
w(t)\\
w_t(t)\\
\bu(t)
\end{bmatrix} \\
[\mathbf{u}(0),w(0),w_t(0)]=[\mathbf{u}_0,w_0,w_1].
\end{cases}
\end{equation}
Here, $\mathbf{A}_\rho : \mathbf{D}(A_\rho)\subset \bH \rightarrow \bH$ generates a strongly continuous semigroup of contractions on $\mathbf{H}$ and is given explicitly as:
\begin{equation}
\label{gener}
\mathbf{A}_{\rho} = \begin{bmatrix}
0 & I & 0\\
-P_\rho^{-1}\Delta^2 + P_\rho^{-1}G_{\rho,1}|_{\Omega_p} & 0 & P_\rho^{-1}G_{\rho,2}|_{\Omega_p}\\
-\nabla G_{\rho,1} & 0 & \nu_f\Delta - \nabla G_{\rho,2}
\end{bmatrix},     
\end{equation}
with its domain $D\left(\mathbf{A}_{\rho}\right)=\left\{\left[w_{1}, w_{2}, \mathbf{u} \right] \subset \mathbf{H}\right.$ defined as

\begin{eqnarray}
w_{1} &\in &\mathcal{S}_{\rho }\equiv \left\{ 
\begin{array}{l}
H^{4}(\Omega _{p})\cap H_{0}^{2}(\Omega _{p})\text{, if }\rho =0 \\ 
H^{3}(\Omega _{p})\cap H_{0}^{2}(\Omega _{p})\text{, if }\rho >0%
\end{array}%
\right. ;  \label{22} \\
w_{2} &\in &H_{0}^{2}(\Omega _{p})\text{, }~~~\mathbf{u}\in \mathbf{H}%
^{2}(\Omega _{f});  \label{22.3} \\
\mathbf{u} &=&\mathbf{0}\text{ \ on }S\text{ and  }\mathbf{u}=[0,0,w_{2}]%
\text{ on }\Omega _{p}.  \label{23}
\end{eqnarray}
(See \cite{A-Clark}.)
Moreover, the associated pressure term $p_0$, say, for $\left[ w_{1},w_{2},\mathbf{u}\right] \in \mathcal{D}(\mathbf{A}_{\rho })$
is given explicitly by the formula in (\ref{pressure}). That is,
\begin{equation}
p_{0}=G_{\rho ,1}(w_{1})+G_{\rho ,2}(\mathbf{u})  \label{p_0} \in H^{1}(\Omega_f).
\end{equation}
(The $H^{1}$-regularity for $p_0$ is inferred from the definition of the $G_i$ in (\ref{G1}) and (\ref{G2}), and the definition of $\mathcal{D}(\mathbf{A}_{\rho })$.)

\section{Main Result: Wellposedness of the Coupled PDE System \eqref{1}-\eqref{4}}\label{sec:wellposedness}
Semigroup generation of the modeling operator $\mathbf{A}_{\rho }:D(\mathbf{A}_{\rho })\subset \mathbf{H}$ was originally established in \cite{A-Clark}. The main goal of this manuscript is to give an alternative and arguably more practical proof (in terms of numerical analysis) to show that the system \eqref{1}-\eqref{4}, or equivalently the abstract ODE system \eqref{ODE}, may be associated with a $C_0-$semigroup $\{e^{\mathbf{A_{\rho}}t}\}_{t\geq 0},$ where $\mathbf{A_{\rho}}:D(\mathbf{A_{\rho}})\subset 
\mathbf{H}\rightarrow \mathbf{H}$ is the matrix operator defined in \eqref{gener}. To this end, we will eventually arrive at an inf-sup system of the classical form \eqref{12.5}, and subsequently put ourselves in position to invoke the Babuska-Brezzi Theorem. We recall the statement of this result for the reader's convenience:
\begin{theorem}
\cite{kesevan}\label{BB} (Babuska-Brezzi) Let $\Sigma ,$ $V$ be real Hilbert spaces and $a:\Sigma \times
\Sigma \rightarrow 
\mathbb{R}
,$ $b:\Sigma \times V\rightarrow 
\mathbb{R}
$ bilinear forms which are continuous. Let%
\[
Z=\left\{ \sigma \in \Sigma :b(\sigma ,v)=0,\text{ \ for every }v\in
V\right\} .
\]%
Assume that $a(\cdot ,\cdot )$ is $Z$-elliptic; i.e., there exists a
constant $\alpha >0$ such that 
\[
a(\sigma ,\sigma )\geq \alpha \left\Vert \sigma \right\Vert _{\Sigma }^{2},%
\text{ \ \ for every }\alpha \in Z.\text{ }
\]%
Assume further that there exists a constant $\beta >0$ such that%
\[
\sup_{\tau \in \Sigma }\frac{b(\tau ,v)}{\left\Vert \tau \right\Vert
_{\Sigma }}\geq \beta \left\Vert v\right\Vert _{V},\text{ \ \ for every }%
v\in V.
\]%
Then if $\kappa \in \Sigma $ and $l\in V,$ there exists a unique pair $%
(\sigma ,u)\in \Sigma \times V$ such that%
\begin{equation}
\left\{ 
\begin{array}{l}
a(\sigma ,\tau )+b(\tau ,u)=(\kappa ,\tau )_{\Sigma} \text{ \ \ \ }\forall \text{ }%
\tau \in \Sigma \\

b(\sigma ,v)=(l,v)_{V}  \text{ \ \ \ \ \ \ \ \ \ \ \ \ \ \ \ \ }\forall \text{ }v\in V.
\end{array}
\right. \label{12.5}
\end{equation}
\end{theorem}

\bigskip

\noindent We begin with the decomposition of the associated pressure variable to
(\ref{1})-(\ref{4}):
\begin{equation}
p(t)=q_{0}(t)+c_{0}(t)  \label{p_qc}
\end{equation}%
where, pointwise in time, $q_{0}\in Q$ and $c_{0}$ is a constant function representing the mean pressure, i.e., $c_0= \int_{\Omega_f} p \; d\Omega_f$.
The mixed variational formulation to be invoked here, quite different from that constructed in \cite{A-Clark}, is such that its solution (or the saddlepoint) is associated with a Lagrange multiplier ${g} \in G := H^{-1/2}(\Omega_p)$ which has the characterization, 
\begin{equation} \label{LMdef}
{g} =  \left. q_0\right\vert _{\Omega _{p}}. 
\end{equation}
That is to say, our particular bilinear form $b(\cdot,\cdot)$ of Theorem \ref{BB}, to be constructed, will involve an accounting of the Dirichlet trace of the FSI pressure on $\Omega_p$. \\

\noindent
We formally derive our mixed variational formulation: \\ Assume that $\left[ w,w_{t},\mathbf{u},p\right] $ is a classical solution of the FSI
system \eqref{1}-\eqref{4}, with $p$ decomposed as in (\ref{p_qc}). Multiplying the governing equations, pointwise in time, by test functions $\left[ \mathbf{v},z\right] \in \mathbf{U}\times H^2_0(\Omega_p)$, and
 subsequently applying integration by parts, we obtain  
\begin{eqnarray}
  \left(  \bu_t, \, \bv \right)_{\Omega_f}
 +  (\nu_{f} \nabla \bu, \nabla \bv)_{\Omega_f} - (q_0,\nabla \cdot
\bv )_{\Omega_f} + \langle g, v_3 \rangle_{\Omega_p} &=&  ( \bbf_f, \bv )_{\Omega_f}   \quad \forall \bv \in  \bU \,\label{st1}, 
\\
\left( w_{tt}, \, z \right)_{\Omega_p} + \rho ( \nabla w_{tt}, \nabla z)_{\Omega_p} + ( \Delta w, \Delta z)_{\Omega_p} -  \langle g+c_0, z\rangle_{\Omega_p} & = & (f_p, z)_{\Omega_p}   \quad \forall z \in  H^2_0(\Omega_p),
\,\label{st2}
\end{eqnarray}
where 
${g}$ is as in (\ref{LMdef}). \\
\noindent
Note that in ariving at this coupled pair of relations we are: (i) using the identification in (\ref{LMdef}); (ii) making the inference that $ \nabla  \cdot \bu =0 $ and $\Omega _{p}$ being contained within the plane $z=0$ means%
\[
\left( \mathbf{\nu }\cdot \nabla \mathbf{u},\mathbf{v}\right) _{\Omega
_{p}}=0\text{, for all }\mathbf{v}\in \mathbf{U}.
\]
Moreover, using the divergence free constraint on the fluid PDE component, and the matching velocities BC in  (\ref{3}) (which we will also come to view here as a constraint), we have
\begin{eqnarray}
 (l, \nabla \cdot \bu)_{\Omega_f}  &=& 0,  \quad \forall l \in Q \label{p1}, 
\end{eqnarray}
where space $Q$ is as defined in (\ref{Q}). In addition, if we set 
\begin{equation}
G\equiv \lbrack H_{\overline{\Omega }_{p}}^{\frac{1}{2}}]^{\prime },
\label{G}
\end{equation}%
the dual of the Hilbert space defined in (\ref{contain}), then using the
matching velocity conditions in (\ref{3}), we then obviously have 
\begin{equation}
\langle w_{t}-u_{3},h \rangle_{\Omega _{p}}=0\text{, \ for every } h \in G.
\label{p2}
\end{equation}
(The reason for the invocation of the space $G$ will be made clear in what follows.) \\

\noindent Collecting the relations (\ref{st1})-(\ref{p1}) and (\ref{p2})-(\ref{p_qc}), and recalling further that the structural displacement is in $\widehat{L}^{2}(\Omega _{p})$, we then see that the above system \eqref{1}-\eqref{3} can be formulated as the following saddle point problem:\\ \\  

\noindent For all $t\in (0,T)$, $\left[ \mathbf{u}(t),w(t),w_t(t)\right] \in \mathbf{U}\times
H_{0}^{2}(\Omega _{p}) \times H_{0}^{2}(\Omega _{p}) $ and $\left[ q_{0}(t),g(t),c_{0}(t)\right] \in 
 Q \times G\times \mathbb{R}$ solve 
\begin{eqnarray}
\label{var1}
A([\bu,w],[\bv,z]) + b([\bv,z], [q_0,g,c_0]) &= & ({ \bbf_f}, \bv )_{\Omega_f} +\left( f_{p},z\right) _{\Omega _{p}},  ~~~ \forall \left[ \mathbf{v},z\right] \in \mathbf{U}\times H_{0}^{2}(\Omega
_{p}), \\ 
b([\bu,w_t], [l, h,r]) & = &0,~~~ \forall \left[ l,h,r\right] \in Q \times G\times \mathbb{R} \label{var2}
\end{eqnarray}
where, for all  $\left[ \mathbf{v},z\right] \in $ $\bU\times H_{0}^{2}(\Omega _{p}),$
$$ A([\bu(t),w(t)],[\bv,z]) :=   \left(  \bu_t, \, \bv \right)_{\Omega_f}
 +  (\nu_{f} \nabla \bu, \nabla \bv)_{\Omega_f} + \left( w_{tt}, \, z \right)_{\Omega_p} + \rho ( \nabla w_{tt}, \nabla z)_{\Omega_p} + ( \Delta w, \Delta z)_{\Omega_p}.$$
Moreover, $b(\cdot,\cdot)$ is defined by having for all  $\left[ \mathbf{U}\times H_{0}^{2}(\Omega _{p})\right] \times \left[ 
 Q \times G\times \mathbb{R}\right] ,$

\begin{equation}
b([\mathbf{v},z],[l,h,r])=-\left( l,\nabla \cdot \mathbf{v}\right) _{\Omega
_{f}}+\langle h,v_{3}\rangle _{\Omega _{p}}-\langle h,z\rangle _{\Omega
_{p}}-\left( r,z\right) _{\Omega _{p}}  \label{b}
\end{equation}
The formally derived variational system (\ref{var1})-(\ref{var2}), will ultimately pave the way for our main result here, which is an alternative and computationally implementable proof of maximality, with respect to the operator $\mathbf{A}_{\rho }:D(\mathbf{A}_{\rho })\subset \mathbf{H}\rightarrow 
\mathbf{H}$, defined in (\ref{gener}). In particular, we have:
 
\begin{theorem}
\label{well}With reference to problem \eqref{1}-\eqref{4} (or Cauchy problem (\ref{ODE})), the operator $\mathbf{A_{\rho}}:D(\mathbf{A_{\rho}})\subset \mathbf{H}%
\rightarrow \mathbf{H}$, defined in (\ref{gener}), generates a $%
C_{0}$-semigroup of contractions on $\mathbf{H}$. Consequently, the solution\\ 
$\Phi (t)=\left[ \bu(t),w(t),w_{t}(t)\right] $ of \eqref{1}-\eqref{4}, or
equivalently (\ref{ODE}), is given by 
\[
\Phi (t)=e^{\mathbf{A_{\rho}}t}\Phi _{0} +\int_{0}^{t}e^{\mathbf{A}_{\rho }(t-s)}\mathbf{F}(s)ds
 \in C([0,T];\mathbf{H})\text{,}
\]%
where $\Phi _{0}=\left[ u_{0},w_{0},w_{1}\right] \in \mathbf{H}$ and $\mathbf{F}=[\mathbf{f}_{f},f_{p}]\in L^{2}((0,T)\times \Omega
_{f})\times L^{2}((0,T)\times \Omega _{p}).$ 
In particular, there exists a mixed variational formulation, devolving from
the formally derived time dependent relations in (\ref{var1})-(\ref{var2}))
which yields the conclusion%
\[
Range(\mathbf{A}_{\rho })=\mathbf{H}.
\]
\end{theorem}

\newpage 

\noindent \textbf{Proof of Theorem \ref{well}}\\

\noindent If one is to invoke the Lumer-Phillips Theorem to justify strongly continuous contraction semigroup generation of  $\mathbf{A_{\rho}}:D(\mathbf{A_{\rho}})\subset \mathbf{H}$, one should show that $\mathbf{A_{\rho}}$ is maximal dissipative. Since the dissipativity of $\mathbf{A}_{\rho }$ (i.e., $\forall \Phi \in D(%
\mathbf{A}_{\rho })$, $\left( \mathbf{A}_{\rho }\Phi ,\Phi \right) _{\mathbf{%
H}}\leq 0)$
has been explicitly shown  in \cite{A-Clark}, we will focus here on the heart of the matter ; i.e., a novel inf-sup approach -- vis-\`{a}-vis
that of \cite{A-Clark} -- to establish the range condition: for $ \lambda > 0,  \left( \lambda \mathbf{I}-\mathbf{A}_{\rho }\right) :D(\mathbf{A}_{\rho
})\rightarrow \mathbf{H}$ is a surjective map.
 \\

\noindent Given  $[\bf_1, f_2, f_3] \in \bH$, we consider the problem of finding $[\bu,w_1,w_2] \in \mathcal{D}(\mathbf{A}_\rho)$ which solves 
\begin{equation}\label{eq4}
    [\lambda\textbf{I}-\mathbf{A}_{\rho}]\begin{bmatrix}
    \bu\\
w_1\\
w_2
\end{bmatrix} = 
\begin{bmatrix}
\bf_1\\
f_2\\
f_3
\end{bmatrix},~~~~~\lambda > 0.
\end{equation}
 The resolvent equation \eqref{eq4} now reads as (taking $\nu_f = 1$), 
 \begin{eqnarray}\label{eq5}
\left\{ \begin{aligned}
&&\lambda \bu- \Delta \bu + \nabla
q_0 = \bbf_1  \quad \mbox{in } \Omega_f\,,  \\
&&\nabla \cdot \; \bu  = 0  \quad \mbox{in } \Omega_f\,, \\
& &\bu = \textbf{0} \quad \mbox{on } S\,, \\
& & \lambda w_1 - w_2 = f_2 \quad \mbox{in } \Omega_p\,, \\
&& \lambda P_{\rho} w_2+\frac{1}{\lambda}\Delta^2 w_2 - q_0|_{\Omega_p} - c_0 = P_{\rho}f_3 - \frac{1}{\lambda}\Delta^2 f_2 \quad \mbox{in } \Omega_p \,,\\
&&w_2 = \frac{\partial w_2}{\partial \bn_p}  = 0  \quad \mbox{on } \partial\Omega_p \,, 
&& \bu = [0, 0, w_2] \quad \mbox{on } \Omega_p,
\end{aligned} \right.
\end{eqnarray}
where $P_{\rho}(.)$ is the operator defined in \eqref{P}. Also, the associated pressure component $p=q_0 + c_0$, where $q_0 \in L^2_0(\Omega_f)$ and $c_0$ is a constant. In order to show the maximality (range) condition for the operator $\mathbf{A_\rho},$ we will use the  static system \eqref{eq5} to derive a variational formulation through  to solve for the variables $[u,p,p|_{\Omega_p},w_2]$; subsequently, we will reconstruct $\{w_1\}$(from $\{w_2\}$ via \eqref{eq5}$_4$). At this point, we note that the abstract resolvent equation \eqref{eq4} is formally a frequency domain version of the time dependent system \eqref{1}-\eqref{4}. In turn, the mixed variational approach which we will eventually formulate will constitute a ``Laplace Transformed'' version of the saddle point problem \eqref{var1}-\eqref{var2}. 
\\ 

\noindent \textbf{The Proof of the Maximality of $\mathbf{A}_\rho$  }
\vspace{0.3cm}

\noindent To start, we multiply the fluid equation \eqref{eq5}$_1$ by $\bv \in \bU$ and the structure equation \eqref{eq5}$_5$ by $z \in H^2_0(\Omega_p)$. Subsequently integrating by parts (and using the fact that normal vector $\left. \nu \right\vert _{\Omega
_{p}}=[0,0,1]$)
 we then obtain,
 \begin{eqnarray}\label{eq7}
\left\{ \begin{aligned}
&&\lambda (\bu, \bv)_{\Omega_f} + (\nabla \bu, \nabla \bv)_{\Omega_f} -
(q_0, \nabla \cdot \bv)_{\Omega_f} + (\left. q_{0}\right\vert _{\Omega _{p}},v_{3})_{\Omega _{p}}= (\bbf_1, \bv)_{\Omega_p},  \\
&& \lambda (w_2, z)_{\Omega_p}+\lambda \rho(\nabla w_2, \nabla z) + \frac{1}{\lambda}(\Delta w_2, \Delta z)_{\Omega_p}
 - (q_0|_{\Omega_p}, z)_{\Omega _{p}} - (c_0, z)_{\Omega_p}\\
 &&= (P_{\rho}f_3, z)_{\Omega_p} - \frac{1}{\lambda}(\Delta f_2, \Delta z)_{\Omega_p}.
\end{aligned} \right.
\end{eqnarray}
If we now set%
\begin{equation}
\Sigma \equiv \mathbf{U}\times H_{0}^{2}(\Omega _{p}),  \label{Sig}
\end{equation}
then the above system \eqref{eq7} gives the following variational formula
\begin{equation}
\left\{ 
\begin{array}{l}
a_{\lambda }([\mathbf{u},w_{2}],[\mathbf{v},z])+b([\mathbf{v}%
,z],[q_{0},g,c_{0}])=F(\mathbf{v},z),~~\text{ \ for }[\mathbf{v},z]\in \Sigma 
\\ 
\\ 
b([\mathbf{u},w_{2}],[l,h,r])=0~~\text{ \ for }Q\times G\times \mathbb{R}%
\end{array}
\right.   \label{eq8}
\end{equation}
Here, the bilinear form $a_{\lambda}(\cdot, \cdot):\bSigma\times \bSigma \rightarrow \real$ is defined by having $\forall $ $[\mathbf{v},z],[\mathbf{\tilde{v}},\tilde{z}]\in \Sigma $,
\begin{equation}
a_{\lambda}([\bv,z],[\tilde{\bv},\tilde{z}]) = \lambda (\bv, \tilde{\bv})_{\Omega_f} + (\nabla v, \nabla \tilde{\bv})_{\Omega_f} + \lambda (z, \tilde{z})_{\Omega_p} + \lambda \rho (\nabla z, \nabla \tilde{z})_{\Omega_p} + \frac{1}{\lambda}(\Delta z, \Delta \tilde{z})_{\Omega_p}.\label{a-l} \end{equation}
Moreover, the (constraint-capturing) bilinear form
$b(\cdot, \cdot): \bSigma \times [Q \times G \times \real] \rightarrow \real$ is given as 
\begin{equation}\label{eq75}
b([\bv, z], [l,h,r]) = -(l, \nabla \cdot \bv)_{\Omega_f}+\langle h,v_3\rangle_{\Omega_p}-\langle h,z\rangle_{\Omega_p} -(r,z)_{\Omega_p},   
\end{equation}
for all $[\mathbf{v},z]\in \Sigma $ and $[l,h,r]\in Q\times G\times \mathbb{R%
}$. In addition, 
\begin{equation}
F(\bv,z)=(\bbf_{1},\bv)_{\Omega _{f}}+(P_{\rho }f_{3},z)_{\Omega _{p}}-\frac{%
1}{\lambda }(\Delta f_{2},\Delta z)_{\Omega _{p}}.  \label{Force}
\end{equation}
It is clear that $a_{\lambda}(\cdot, \cdot)$ is continuous and $\bSigma-$elliptic. Hence, our main task here is to show that $b(\cdot, \cdot)$ satisfies the Babuska-Brezzi inequality, with respect to the given choice of spaces. That is,  we must verify there exists $ \beta > 0$ such that $\forall~ [l,h,r] \in Q\times G\times \mathbb{R}$,
\begin{equation}\label{eq9}
\sup_{[\bv,z] \in \bSigma}\frac{b([\bv, z], [l,h,r])}{||[\bv,z]||_{\bSigma}} \geq \beta ||[l,h,r]||_{Q\times G\times \mathbb{R}}.    
\end{equation}
We establish the inequality (\ref{eq9}) through several steps:
\paragraph{Step I:}
Let $[l,h,r]\in Q\times G\times \mathbb{R}$
  be given. Therewith, let $\bmu_0 \in H^1_{0}(\Omega_f)$ denote the solution of 
\begin{eqnarray}\label{eq10}
\left\{ \begin{aligned}
&& div(\bmu_0) = -l \quad \mbox{in } \Omega_f, \\
&& \bmu_0|_{\partial \Omega_f} = 0 \quad \mbox{on } \partial \Omega_f. 
\end{aligned} \right.
\end{eqnarray}
Since $l \in Q$,  it is well known (see e.g.,,  Lemma 3.2.2 of \cite{kesevan}) that $\exists \bmu_0 \in H^1(\Omega_f)$ which solves \eqref{eq10} and obeys the estimate 
\begin{equation}\label{eq11}
    ||\nabla \bmu_0||_{\Omega_f} \leq C_0 ||l||_{\Omega_f}.
\end{equation}
\paragraph{Step II:} By Hahn-Banach Theorem, for a given data $h \in G$,  there exists a function $k \in H^{1/2}_{\overline{\Omega_p}}$ such that
\begin{eqnarray}\label{eq12}
\left\{ \begin{aligned}
&& \langle h, k \rangle_{\Omega_p} = ||h||_{G}, \\
&& ||k||_{H^{1/2}_{\overline{\Omega_p}}} = 1 
\end{aligned} \right.
\end{eqnarray}
(see Corollary 6.2.10 in \cite{hutson}). We note super-critically that $g$ can be extended by zero to all of $\partial \Omega_f$ so as to attain a $H^{1/2}(\partial \Omega_f)$ function [see Theorem 3.29 (ii) and Theorem 3.33 of \cite{McLean}]. Consequently, the following BVP in solution variable $\bmu_1 \in \bU$ is wellposed: 
\begin{eqnarray}\label{eq13}
\left\{ \begin{aligned}
&& div(\bmu_1) = \frac{\eta}{meas(\Omega_f)} ||h||_{G} \int_{\Omega_p}k~\partial \Omega_p \quad \mbox{in } \Omega_f \\
&& \bmu_1 = \begin{cases} 
      0 & \mbox{on } S \\
      [0,0,\eta \left(||h||_{G}\right)k] & \mbox{on } \Omega_p, 
   \end{cases}
\end{aligned} \right.
\end{eqnarray}
where $\eta > 0$ is a parameter to be determined. (Note that the compatibility condition is satisfied for the solvability of the problem \eqref{eq13}.) Moreover, by \eqref{eq12}, we have 
\begin{equation}\label{eq14}
    ||\nabla \bmu_1||_{\Omega_f} \leq C_2||h||_{G}.
\end{equation}
(Here, we neglect noting the dependence of constant $C_2$ on given $\eta$, as it is ultimately benign; see \eqref{eq20}.) 
\paragraph{Step III:} We define the realization of the biharmonic operator on $\Omega_p$ with clamped boundary conditions: 
 $$\textbf{\AA }
:\mathcal{D}(\textbf{\AA }
) \subset L^2(\Omega_p)\rightarrow L^2(\Omega_p);~~~\textbf{\AA }
\psi = \Delta^2 \psi, ~~~\forall \psi \in \mathcal{D}(\textbf{\AA }
),$$ 
 where $\mathcal{D}(\textbf{\AA }
)=H^4(\Omega_p)\cap H_0^2(\Omega_p)$. 
 \begin{remark}
Note that we are implicitly invoking elliptic theory in specifying the domain $\mathcal{D}(\textbf{\AA }
)$ if $\Omega_p$ has smooth boundary, or Theorem 2,  p. 563, of \cite{BR} if $\Omega_p$ is a convex polygon with angles $\theta \leq \frac{2\pi}{3}$.
 \end{remark}
\noindent Now, since $\textbf{\AA }$ is a self-adjoint, positive definite operator we have by \cite{Grisvard} the characterization
\begin{equation} 
    \mathcal{D}(\textbf{\AA }
^{\frac{1}{2}}) \approx H^2_0(\Omega). \label{charA}
\end{equation}
Moreover, the dual space $H^{-1/2}(\Omega_p)$ is an isometric realization of $[H^{1/2}_{\overline{\Omega_p}}]' (=G)$ (see Theorem 3.29 (ii) and Theorem 3.14 (ii) of \cite{McLean}). Therefore, there exists $ C^{*}>0$ such that
\begin{equation}\label{eq15}
    ||\textbf{\AA }
^{- \frac{1}{2}}f||_{\Omega_p} = ||f||_{H^{-2}(\Omega_p)} \leq C^{*} ||f||_{G}, \quad \quad \quad \forall ~f \in G.
\end{equation}
Here, we are using the fact from (\ref{charA}) that $$||f||_{H^{-2}(\Omega_p)} = ||\textbf{\AA }
^{-\frac{1}{2}}f||_{\Omega_p}, \quad \quad \forall f \in H^{-2}(\Omega_p).$$
\paragraph{Step IV:} We set now
\begin{eqnarray}
\bmu &=&\bmu_{0}+\bmu_{1};  \label{eq16} \\
\phi  &=&-\text{\textbf{\AA }}^{-1}(h+r),  \nonumber
\end{eqnarray}
where the $\bmu_i$ are from (\ref{eq10}) and (\ref{eq13}). Then, from \eqref{eq15} and \eqref{eq16}, we have
\begin{equation}\label{eq17}
  ||\textbf{\AA }^{-\frac{1}{2}}
\phi||_{\Omega_p} \leq C\left[||h||_{G} + |r|\right]. 
\end{equation}
Also, taking into account \eqref{eq75}, \eqref{eq10} and \eqref{eq13} we get 
\begin{eqnarray} \label{eq18}
b([\bmu, \phi], [l,h,r]) &&= -(l, div(\bmu_0))_{\Omega_f} + \langle h,(\bmu|_{\Omega_p})_{3}\rangle_{\Omega_p} + \langle h+r,\text{\textbf{\AA }}^{-1}(h+r)\rangle_{\Omega_p}, \nn \\
&&(\text{after using the fact that } l \in Q  )  \nn \\ 
&&= ||l||^2_{\Omega_f} + \eta ||h||_{G}+||\textbf{\AA }^{-\frac{1}{2}}(h+r)||^2_{\Omega_p} \nn\\
&& = ||l||^2_{\Omega_f} + \eta ||h||_{G}+||\textbf{\AA }^{-\frac{1}{2}}h||^2_{\Omega_p} \nn\\
&& + \boxed{2\left( \text{\textbf{\AA }}^{-\frac{1}{2}}h,\text{\textbf{\AA }}^{-\frac{1%
}{2}}r\right) _{\Omega _{p}}}+||\textbf{\AA }^{-\frac{1}{2}}r||^2_{\Omega_p}
\end{eqnarray}
Now, via 
\begin{equation*}
    2|ab| \leq \epsilon a^2 + \frac{b^2}{\epsilon}, \quad 0< \epsilon < 1,
\end{equation*}
we have for the boxed term in \eqref{eq18},
\begin{equation}
\label{*}
    2\left( \text{\textbf{\AA }}^{-\frac{1}{2}}h,\text{\textbf{\AA }}^{-\frac{1%
}{2}}r\right) _{\Omega _{p}} \geq -\frac{\left\Vert \text{\textbf{\AA }}^{-\frac{1}{2}}h\right\Vert _{\Omega
_{p}}^{2}}{\epsilon }-\epsilon \left\Vert \text{\textbf{\AA }}^{-\frac{1}{2}%
}r\right\Vert _{\Omega _{p}}^{2}
.
\end{equation}
So, for given $[l,h,r] \in Q \times G\times \real$, upon incorporating \eqref{*} into \eqref{eq18}, we arrive at 
\begin{eqnarray} \label{crux}
\sup_{[\bv,z] \in \bSigma}\frac{b([\bv, z], [l,h,r])}{||[\bv,z]||_{\bSigma}} &&\geq \frac{b([\bmu, \phi], [l,h,r])}{||\nabla \bmu||_{\Omega_f} +  |\left\Vert \text{\textbf{\AA }}^{\frac{1}{2}}\phi \right\Vert _{\Omega _{p}} }\\
&& \geq \frac{||l||^2_{\Omega_f} + \eta ||h||^2_{G}+(1-\epsilon)||\text{\textbf{\AA }}^{\frac{1}{2}}r||^2_{\Omega_p} - \left(\frac{1}{\epsilon}-1\right)||\text{\textbf{\AA }}^{\frac{1}{2}}h||^2_{\Omega_p}}{||\nabla \bmu||_{\Omega_f} +  |\left\Vert \text{\textbf{\AA }}^{\frac{1}{2}}\phi \right\Vert _{\Omega _{p}} }\\
&& \geq \frac{||l||^2_{\Omega_f} + \left(\eta -C^{*}(\frac{1}{\epsilon }-1) \right) ||h||^2_{G}+(1-\epsilon)||\text{\textbf{\AA }}^{\frac{1}{2}}r||^2_{\Omega_p}}{||\nabla \bmu||_{\Omega_f} + \left\Vert \text{\textbf{\AA }}^{\frac{1}{2}}\phi \right\Vert _{\Omega _{p}}  },
\end{eqnarray}
where in obtaining this intermediate inequality, we also used \eqref{eq15}. If we choose $\eta > C^{*}(\frac{1}{\epsilon }-1)$ we then obtain
\begin{eqnarray}\label{eq20}
\sup_{[\bv,z] \in \bSigma}\frac{b([\bv, z], [l,h,r])}{||[\bv,z]||_{\bSigma}} 
\geq \frac{C_q \left(||l||^2_{\Omega_f} + ||h||^2_{G}+|r|^2 \right)}{||\nabla \bmu||_{\Omega_f} +  ||\text{\textbf{\AA }}^{\frac{1}{2}}\phi||_{\Omega_p}}
\end{eqnarray}
where $$C_q = max\{1, \left(\eta -C^{*}(\frac{1}{\epsilon }-1)\right), (1-\epsilon)||\text{\textbf{\AA }}^{\frac{1}{2}}(1)||^2_{\Omega_p}\}.$$ Continuing, with the estimates \eqref{eq11}, \eqref{eq14} and \eqref{eq17} in mind, we get 
\begin{eqnarray}\label{eq21}
\sup_{[\bv,z] \in \bSigma}\frac{b([\bv, z], [l,h,r])}{||[\bv,z]||_{\bSigma}} 
&& \geq \frac{\tilde{C} \left(||l||^2_{\Omega_f} + ||h||^2_{G}+|r|^2 \right)}{||l||_{\Omega_f} + ||h||_{G}+|r|} \nn \\
&& \geq \left(\frac{\tilde{C}}{3}\right) \frac{ \left(||l||_{\Omega_f} + ||h||_{G}+|r| \right)^2}{||l||_{\Omega_f} + ||h||_{G}+|r|} \nn \\
&&\geq \left(\frac{\tilde{C}}{3}\right)\left(||l||_{\Omega_f} + ||h||_{G}+|r| \right).
\end{eqnarray}
This last estimate gives the conclusion that the Babuska-Brezzi estimate \eqref{eq9} is satisfied with $\beta = \frac{\tilde{C}}{3}$. 
In sum: with respect to the inf-sup system \eqref{eq8}, one has wellposedness of
the solution map%
\begin{equation}
\left[ \mathbf{f}_{1},f_{2},f_{3}\right] \in \mathbf{H}\rightarrow \left\{ \lbrack 
\mathbf{u},w_2],[q_{0},g,c_{0}]\right\} \in \lbrack \mathbf{U}\times
H_{0}^{2}(\Omega _{p})]\times \lbrack Q\times G\times \mathbb{R}].
\label{star}
\end{equation} 

\newpage 

\noindent\textbf{Completion of the Proof of Theorem \ref{well}}
\vspace{0.3cm}

\noindent To conclude the proof of Theorem \ref{well}, we will show that the solution  $[\bu,w_2] \in \Sigma$ of \eqref{eq8} gives rise indeed to a solution $[\bu,w_1, w_2]$ to the static PDE system \eqref{eq5}, and which is in the domain of the FSI generator $D(\mathbf{A}_\rho)$ in \eqref{gener}. To start, we write out explicitly the variational relation in \eqref{eq8}: 
\begin{eqnarray}  \label{eq23}
\left\{ \begin{aligned}
&&\lambda (\bu, \bv)_{\Omega_f} +(\nabla \bu, \nabla \bv)_{\Omega_f} -
(q_0, \nabla \cdot \bv)_{\Omega_f} + \langle g, v_3\rangle_{\Omega_p}  + \lambda (w_2, z)_{\Omega_p} \\
&& +\lambda \rho(\nabla w_2, \nabla z) + \frac{1}{\lambda}(\Delta w_2, \Delta z)_{\Omega_p}
 - \langle g, z\rangle_{\Omega_p}   - (c_0, z)_{\Omega_p}\\
&&= (\bbf_1, \bv)_{\Omega_f}   + (P_{\rho}f_3, z)_{\Omega_p} - \frac{1}{\lambda}(\Delta f_2, \Delta z)_{\Omega_p}, \quad \forall [\bv, z] \in \bSigma.
\end{aligned} \right.
\end{eqnarray}
If we also invoke the constraint equation in \eqref{eq8} (where $b(\cdot,\cdot)$ is as defined in \eqref{eq75}), and choose $h=r=0$ therein, we have 
\begin{equation}\label{eq23a}
    ( div (\bu), l )_{\Omega_f} = 0, \quad \forall l \in Q.
\end{equation}
 Also, taking $l=0$ and $r=0$ in \eqref{eq8} we get 
\begin{equation}\label{eq24}
    u_3 = w_2 \quad \mbox{in } \Omega_p.
\end{equation}
 Lastly, taking $l=h=0$ in \eqref{eq8} gives
 \begin{equation}\label{eq25}
     \int_{\Omega_p}w_2 d\Omega_p = 0.
 \end{equation}
In turn, \eqref{eq24} and \eqref{eq25} yield,
\begin{equation}\label{eq26}
  \int_{\Omega_p}u_3 d\Omega_p =   \int_{\Omega_p}w_2 d\Omega_p = 0.
\end{equation}
Consequently, by Green's Theorem we obtain 
\begin{equation}\label{eq27}
 \int_{\Omega_p}div(\bu) d\Omega_f = \int_{\partial \Omega_f} \bu \cdot \nu d\Omega_f = \int_{\Omega_p}u_3 d\Omega_p = 0  
\end{equation}
Then \eqref{eq27} and \eqref{eq23a} yield that 
\begin{equation}\label{eq28}
    div(\bu) = 0. 
\end{equation}
 Moreover, taking $\bv \in [\mathcal{D}(\Omega_f)]^3$ and $ z = 0 $ in \eqref{eq23} gives
\begin{eqnarray}
\begin{aligned}
\lambda (\bu, \bv)_{\Omega_f}+ (\Delta \bu,  \bv)_{\Omega_f} +
(\nabla q_0, \bv)_{\Omega_f} = (\bbf_1, \bv)_{\Omega_f}, 
\end{aligned}.
\end{eqnarray}
which implies (after the usual density argument) that
\begin{equation}\label{eq29}
    \lambda \bu - \Delta \bu + \nabla q_0 = \bf_{1}.
\end{equation}
Subsequently, using the surjectivity of the bounded Dirichlet trace map $\mathbf{H}^1(\Omega_f) \rightarrow \mathbf{H}^{\frac{1}{2}}(\partial \Omega _{f})$, the wellposedness result (\ref{star}), and Green's Theorem, we have the estimate
\begin{equation}\label{eq30}
||\frac{\partial \bu}{\partial \nu} - q_0\nu||_{H^{-1/2}(\partial \Omega_f)} \leq C||[\bbf_1,f_2,f_3]||_{\bH}.
\end{equation}
In turn, we  take $\bv \in \bU$ and $z=0$ in \eqref{eq23}; this gives,
\[
\lambda (\mathbf{u},\mathbf{v)}_{\Omega _{f}}+(\nabla \mathbf{u},\nabla 
\mathbf{v)}_{\Omega _{f}}-(q_{0},\nabla \cdot \mathbf{v})_{\Omega
_{f}}+ \langle g,v_{3} \rangle_{\Omega _{f}}=(\mathbf{f}_{1},\mathbf{v})_{\Omega _{f}.}.
\]
Integrating by parts, with respect to this relation and invoking \eqref{eq29}-\eqref{eq30}, we have now,
\begin{equation}
    \langle g-q_0, v_3 \rangle_{\Omega_p} = 0, \quad \forall \bv \in \bU. \label{gq}
\end{equation}
With this relation in hand, an appeal to \cite[Theorem 3.29 (ii) and Theorem 3.33]{McLean}, and use of the fact that the Sobolev Trace mapping $\gamma_0: H^1(\Omega_f)\rightarrow H^{1/2}(\partial \Omega_f)$ is surjective, and so has a continuous right inverse $\gamma_0^+(.)$, give the conclusion that
\begin{equation}\label{eq31}
    \left. q_{0}\right\vert _{\Omega _{p}}=g.
\end{equation}
\underline{Indeed:} if projection $\mathbb{P}_3\bv = v_3$ for given $\bv \in \bU$, and if $\gamma_0 \in \mathcal{L}(H^1(\Omega_f), H^{1/2}(\partial\Omega_f))$ satisfies 
\begin{equation*}
    \gamma_0(f) = f|_{\partial \Omega_f}, \quad \forall f \in H^1(\Omega_f),
\end{equation*}
then $[\gamma_0(\mathbb{P}_3(\bU))] \subset H^{1/2}_{\overline{\Omega_p}}.$
On the other hand , if $z \in H^{1/2}_{\overline{\Omega_p}}$, then by Theorem 3.33 of \cite{McLean} it can be extended by zero onto all of $H^{1/2}(\partial \Omega_f)$. Thus if $\gamma_0^{+} \in \mathcal{L}(H^{1/2}(\partial\Omega_f),H^1(\Omega_f) )$ is the continuous right inverse of $\gamma_0$, then $\gamma_0^{+}(z_{ext}) \in H^1(\Omega_f)$ where
\begin{equation}
z_{ext} = \begin{cases}
z \quad \mbox{on } \Omega_p \\
0 \quad \mbox{on } \partial \Omega_f / \Omega_p 
\end{cases} 
\end{equation}
Thus, 
\[
z=\gamma _{0}(\gamma _{0}^{+}(z_{ext}))=\gamma _{0}\left( \mathbb{P}_{3}%
\left[ 
\begin{array}{c}
0 \\ 
0 \\ 
\gamma _{0}^{+}(z_{ext})%
\end{array}%
\right] \right) ,
\]
and so $[\gamma_0(\mathbb{P}_3(\bU)]_{\Omega_p} = H^{1/2}_{\overline{\Omega_p}}$. This characterization and \eqref{gq} give \eqref{eq31}.\\

\noindent For the mechanical component of \eqref{eq23}: Via \eqref{eq31}, the relation \eqref{eq23} with $\bv=0$ gives
\begin{eqnarray*}
\begin{aligned}
&& \lambda (w_2, z)_{\Omega_p}+\lambda \rho(\nabla w_2, \nabla z) + \frac{1}{\lambda}(\Delta w_2, \Delta z)_{\Omega_p}
 - ([q_0|_{\Omega_p} + c_0], z)_{\Omega_p}\\
&&= (P_{\rho}f_3, z)_{\Omega_p} - \frac{1}{\lambda}(\Delta f_2, \Delta z)_{\Omega_p} \quad \forall z \in H^2_0(\Omega_p).
\end{aligned}
\end{eqnarray*}
An integration by parts and density argument then yield,
\begin{equation}
    \lambda w_2 - \lambda \rho \Delta w_2 + \frac{1}{\lambda} \Delta^2 w_2 -[q_0|_{\Omega_p} + c_0]=  P_{\rho}f_3 - \frac{1}{\lambda}\Delta^2 f_2. 
\end{equation}
Also, via the resolvent relation \eqref{eq4}, 
\begin{equation}\label{resolve}
  \lambda w_1 = w_2 + f_2,
\end{equation}
and so we have
\begin{equation}\label{eq34}
  \lambda P_{\rho}w_2 + \Delta^2 w_1 - p|_{\Omega_p} = P_{\rho}f_3. 
\end{equation}
Here again,
\begin{equation}\label{eq35}
 p = q_0 + c_0.  
\end{equation}
(And so the aforementioned Lagrange multiplier in \eqref{LMdef} is $\mathfrak{g}=\left. q_{0}\right\vert _{\Omega _{p}} + c_0$.)

\medskip

\noindent Finally, collecting \eqref{star}, \eqref{eq24}, \eqref{eq25}, \eqref{eq28}, \eqref{eq29}, \eqref{resolve}, \eqref{eq34}, and moreover, using the (geometric) Assumption \ref{geometry}, we can invoke the higher regularity results established in \cite{avalos2014exponential} -- see Lemma 3.1 and Appendix, p. 74 therein -- to conclude that the solution variables $[\bu, w_1,w_2, p]$ are in $\mathcal{D}(A_\rho)$ as given in \cite[page 7]{A-Clark}. In particular, we have that the solution of \eqref{eq5} or \eqref{eq8} satisfies $\bu \in \bH^2(\Omega_f), p \in H^1(\Omega_f)$,
\begin{eqnarray}
w_1 = \left\{ \begin{aligned}
&& H^3(\Omega_p)\cap H^2_0(\Omega_p),  \quad \mbox{if } \rho > 0 \\
&& H^4(\Omega_p)\cap H^2_0(\Omega_p),  \quad \mbox{if } \rho = 0 
\end{aligned} \right.
\end{eqnarray}
and $w_2 \in H^2_0(\Omega_p)$. In short, weak solution $[\bu, w_1,w_2, p]$ of \eqref{eq23} is a classical solution of system \eqref{eq5}. This finishes the proof of Theorem \ref{well}.

 \medskip

\section{Numerical results}\label{sec:NR}
In this section, we consider a benchmark problem from \cite{A-Clark} that addresses fluid-structure interaction in the case of $\rho = 0$. The fluid domain is defined as $\Omega_f:= [0,1] \times[0,1] \times [-1,0]$, and the plate domain $\Omega_p := [0,1] \times [0,1] \times \{0\}$ lies along the top boundary of the fluid domain in the $x y$ plane, as illustrated in \Cref{fig:3d2d}. 
\newcommand{\Depth}{8}
\newcommand{\Height}{4}
\newcommand{\Width}{4}
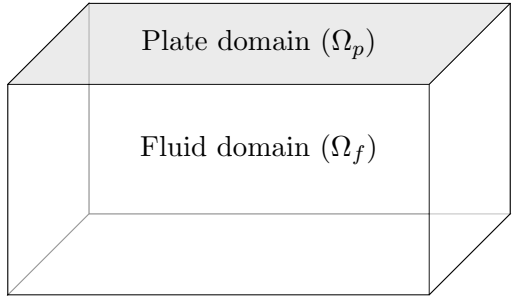
\begin{figure}
\begin{center}
\begin{tikzpicture}[scale=0.7]
\centering
\coordinate (O) at (0,0,0);
\coordinate (A) at (0,\Width,0);
\coordinate (B) at (0,\Width,\Height);
\coordinate (C) at (0,0,\Height);
\coordinate (D) at (\Depth,0,0);
\coordinate (E) at (\Depth,\Width,0);
\coordinate (F) at (\Depth,\Width,\Height);
\coordinate (G) at (\Depth,0,\Height);
\coordinate (H) at (4,2,2);
\coordinate (I) at (4,4,2);

\draw[line width=0.01mm, black,fill=white!20] (O) -- (C) -- (G) -- (D) -- cycle;
\draw[line width=0.01mm,black,fill=white!30] (O) -- (A) -- (E) -- (D) -- cycle;
\draw[line width=0.01mm,black,fill=white!10] (O) -- (A) -- (B) -- (C) -- cycle;
\draw[line width=0.01mm,black,fill=white!20,opacity=0.8] (D) -- (E) -- (F) -- (G) -- cycle;
\draw[line width=0.01mm,black,fill=white!20,opacity=0.6] (C) -- (B) -- (F) -- (G) -- cycle;
\draw[line width=0mm,black,fill=gray!20,opacity=0.8] (A) -- (B) -- (F) -- (E) -- cycle;

    \node at (H) {Fluid domain $(\Omega_f)$};
    \node at (I) {Plate domain $(\Omega_p)$};
\end{tikzpicture}
\caption{3D Fluid and 2D plate system}
\label{fig:3d2d}
\end{center}
\end{figure}
The numerical example in \cite{A-Clark} is for the steady-state case. Accordingly, the fluid problem is approximated using the time-discretized equation:
$$\lambda \bu - \Delta \bu + \nabla p= \bf_{1} $$
where $\lambda$ is the reciprocal of the time step $\Delta t$.
For the plate equation \eqref{2} and the coupling condition \eqref{3}, we set $w_1 :=w$ and 
introduce an additiomal variable $w_2$ representing $w_t$. The time discretized form of \eqref{2} then becomes
\begin{align*}
\lambda w_2 + \Delta^2 w_1 &= p|_{\Omega_p},\\
\lambda w_1-w_2 &= f_2.
\end{align*}
As in \cite{A-Clark}, the following manufactured solution is employed with $\lambda=1$
$$
\begin{aligned}
w_{1}^{\mbox{exact}}= & -x^{4}(x-1)^{4}(2 x-1) y^{4}(y-1)^{4} \\
w_{2}^{\mbox{exact}}= & -\Delta w_{1}^{\mbox{exact}}=12 x^{2}(x-1)^{2}(2 x-1)\left(6 x^{2}-6 x+1\right) y^{4}(y-1)^{4}+ \\
& \quad x^{4}(x-1)^{4}(2 x-1) 4 y^{2}(y-1)^{2}\left(14 y^{2}-14 y+3\right) \\
u_{1}^{\mbox{exact}}= & {\left[2 x^{3}(x-1)^{3}\left(9 x^{2}-9 x+2\right) y^{4}(y-1)^{4}+\right.} \\
& \left.\quad(4 / 5) x^{5}(x-1)^{5} y^{2}(y-1)^{2}\left(14 y^{2}-14 y+3\right)\right]\left[-30 z^{4}-60 z^{3}-30 z^{2}\right] \\
u_{2}^{\mbox{exact}}= & 0 ; \\
u_{3}^{\mbox{exact}}= & -\left[12 x^{2}(x-1)^{2}(2 x-1)\left(6 x^{2}-6 x+1\right) y^{4}(y-1)^{4}+\right. \\
& \left.\quad 4 x^{4}(x-1)^{4}(2 x-1) y^{2}(y-1)^{2}\left(14 y^{2}-14 y+3\right)\right]\left[-6 z^{5}-15 z^{4}-10 z^{3}-1\right] \\
p^{\mbox{exact}}= & 0,
\end{aligned}
$$
and the corresponding right-hand side functions are derived from this solution:
$$
\begin{aligned}
\bf_{1} & = \lambda u^{\mbox{exact}}-\Delta u^{\mbox{exact}} + \nabla p^{\mbox{exact}}, 
f_{2} = \lambda w_{1}^{\mbox{exact}}-w_{2}^{\mbox{exact}}.
\end{aligned}
$$
Notice that the velocity $u=\left[u_{1}, u_{2}, u_{3}\right]$ is divergence free and satisfies $u=[0,0,0]$ and $\Delta u \cdot \nu=[0,0,0]$ on $S$. In the structure domain $\Omega_p$, we have $u=\left[0,0, w_{2}\right]$. 
Furthermore, the pressure $p=0$ 
because $w_2$ is defined as  $-\Delta w_{1}$, which leads to cancellation of the  two terms. Additionally, 
$w_{1}$ and $u$ are chosen such that $\Delta^{2} w_{1}=-\Delta u \cdot \nu$ on $\Omega_p$, and $\Delta u \cdot \nu=0$ on $S$. 

\begin{figure}[!h] 
    \centering
    \begin{subfigure}[t]{0.5\textwidth}
        \centering
    \includegraphics[scale=0.3]{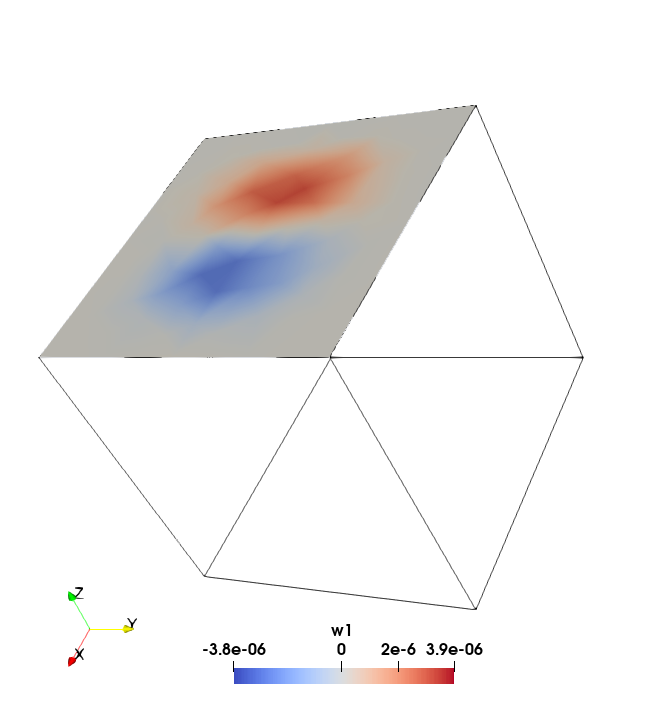} 
        \caption{$w_1$}
    \end{subfigure}%
    \hspace{-0.5in} 
    \begin{subfigure}[t]{0.5\textwidth}
        \centering
      \includegraphics[scale=0.3]{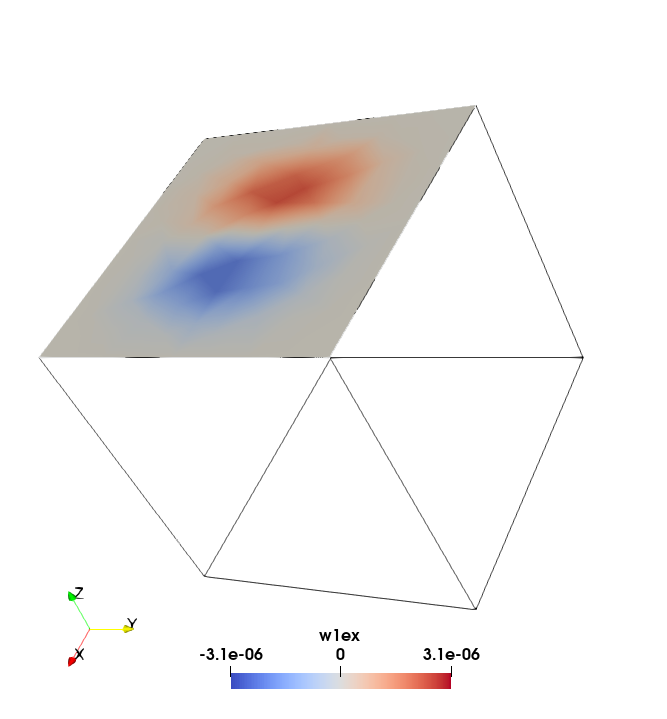}
        \caption{$w_1^{\mbox{exact}}$}
    \end{subfigure}
\caption{$w_1$-FEM (left) and $w_1$-Exact (right) for $h = 1/12$.}
 \label{fig:w1}
\end{figure}

\begin{figure}[!h] 
    \begin{subfigure}[t]{0.5\textwidth}
        \centering
      \includegraphics[scale=0.3]{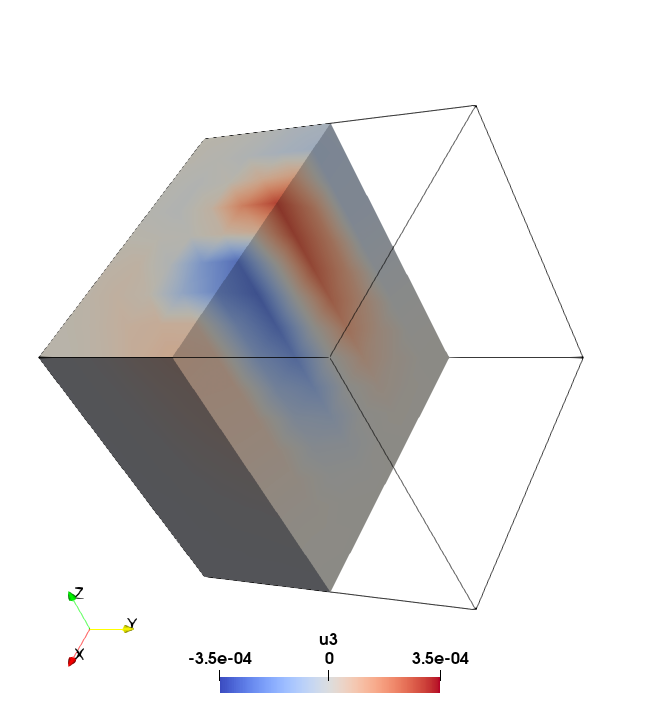}
        \caption{$u_3$}
    \end{subfigure}
    \hspace{-0.5in}
    \begin{subfigure}[t]{0.5\textwidth}
        \centering
      \includegraphics[scale=0.3]{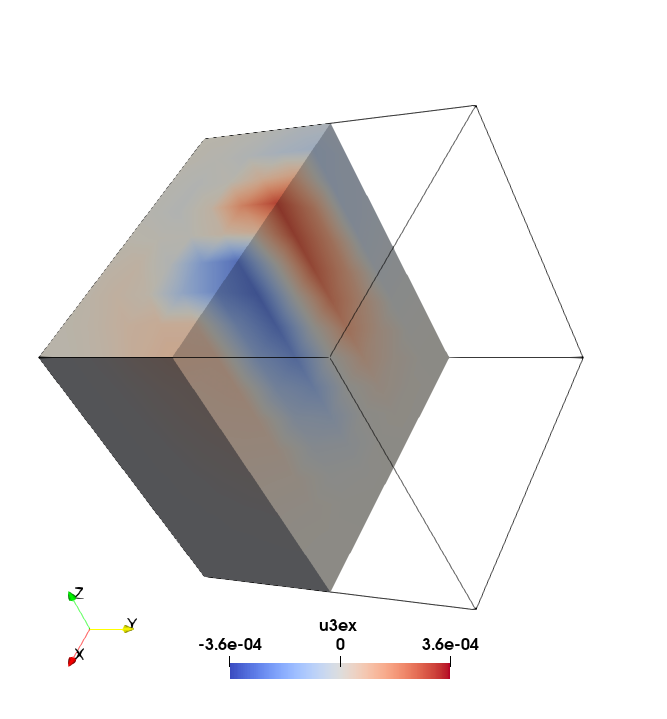}
        \caption{$u_3^{\mbox{exact}}$}
    \end{subfigure}
\caption{$u_3$-FEM (left) and $u_3$-Exact (right)  for $h = 1/12$.}
 \label{fig:u3}
\end{figure}

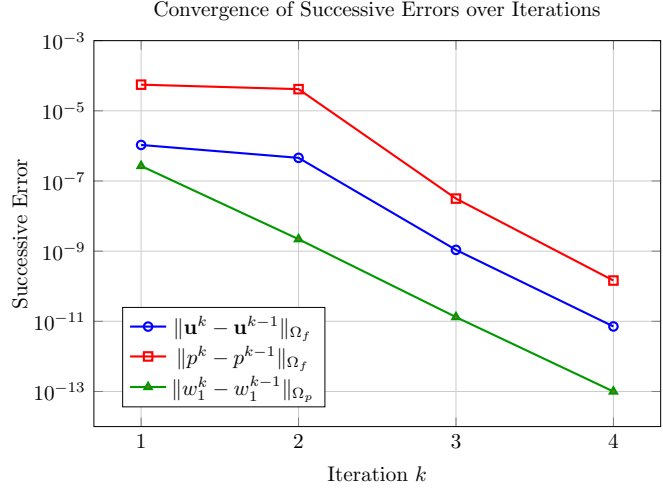
\begin{figure}[h!]
\centering
\begin{tikzpicture}[scale=0.8]
\begin{semilogyaxis}[
    width=11cm,
    height=8cm,
    xlabel={Iteration $k$},
    ylabel={Successive Error},
    title={Convergence of Successive Errors over Iterations},
    legend style={at={(0.05,0.05)}, anchor=south west},
    xtick={1,2,3,4},
    ymin=1e-14, ymax=1e-3,
    grid=both,
    major grid style={line width=.2pt,draw=gray!40},
    minor grid style={line width=.1pt,draw=gray!20},
    mark options={solid},
    every axis plot/.append style={line width=1pt},
    font=\small,
]

\addplot+[mark=o, blue] coordinates {
    (1,1.06353e-06)
    (2,4.56211e-07)
    (3,1.08685e-09)
    (4,7.13518e-12)
};
\addlegendentry{$\|\bu^k - \bu^{k-1}\|_{\Omega_f}$}

\addplot+[mark=square, red] coordinates {
    (1,5.55916e-05)
    (2,4.14472e-05)
    (3,3.14985e-08)
    (4,1.45357e-10)
};
\addlegendentry{$\|p^k - p^{k-1}\|_{\Omega_f}$}

\addplot+[mark=triangle*, green!60!black] coordinates {
    (1,2.6759e-07)
    (2,2.19277e-09)
    (3,1.31924e-11)
    (4,9.9913e-14)
};
\addlegendentry{$\|w_1^k - w_1^{k-1}\|_{\Omega_p}$}

\end{semilogyaxis}
\end{tikzpicture}
\caption{Convergence of successive errors for $\bu$, $p$, and $w_1$ over iterations.}
\label{fig:successive_errors}
\end{figure}

For the finite elements simulation,  we used $(\textbf{P}_2, P_1, P_{2}\mbox{Morley},P_2)$ elements for $(\bu, p, w_1, w_2)$. We evaluated errors and convergence rates for different mesh size $(h)$ (see Table \ref{3dtab:2nd} and Table \ref{3dtab:2ndw}). The results indicate that the $H^1$- and $H^2$- rates of the structure are consistent with the theoretical rates established in \cite{Morleyrate}, while $L^2$ rates are suboptimal for smaller $h$. 
Figures \ref{fig:w1} and \ref{fig:u3} provide a comparison between the finite element solutions ($w_1$-FEM and $u_3$-FEM)  and the corresponding exact solutions ($w_1$-Exact and $u_3$-Exact).
The close agreement observed between the numerical and exact solutions, with both $w_1$ and $u_3$ exhibiting similar behavior, highlights the accuracy of the finite element approximation and the consistency of the numerical model with the theoretical predictions. Further, to assess the performance of the proposed algorithm, we examine the successive errors in the computed solutions for the velocity $\bu$, pressure $p$, and plate functions $w_1$ and $w_2$. These errors are measured at each iteration and provide insight into the convergence behavior of the method. Our results demonstrate that the errors in all variables consistently decrease as the iteration count increases (see \Cref{fig:successive_errors}) confirming the effectiveness of the iterative scheme in accurately resolving the coupled system.

\begin{algorithm}[hbt!]
\SetAlgoLined
\textbf{Input:} $\dot{w}$ initial guess, $\epsilon$ tolerance and $N_{\mbox{iter}}$ maximum number of iteration.\\
\textbf{Output:} $w_2^{k}$\\
$k=0$, error $> \epsilon$, $w_2^{0}=0, \bu^{0}=\textbf{0}$ \\
\textbf{while} $k<N_{\mbox{iter}}$ and error $> \epsilon,$ \textbf{do}

\begin{enumerate}
    \item Solve for $(\bu^{k},p^{k}) \in \bU \times Q$ with $u_3^{k} = w_2^{k-1}$:
\begin{eqnarray}
\begin{aligned}
\lambda \bu^{k} - \Delta \bu^{k} + \nabla p^{k} &= \bf_{1}^k  \quad \mbox{in } \Omega_f \,, 
\\
\nabla \cdot \; \bu^{k} & = 0  \quad \mbox{in } \Omega_f \,, \\
\bu^{k} & = \textbf{0} \quad \mbox{on } S, \\
\bu^{k} = [u_1^{k}, u_2^{k}, u_3^{k}] &= [0, 0, w_2^{k-1}] \quad \mbox{on } \Omega_p,\\
\int_{\Omega_f} p^{k} \; d \Omega_f &=0.
\end{aligned}
\end{eqnarray}    
\item Solve for $(w_1^{k}, w_2^{k})$ in $H^2_0(\Omega_p) \times H^2_0(\Omega_p)$:
 \begin{eqnarray}
\begin{aligned}
\lambda w_2^{k} + \Delta^2 w_1^{k} &= p^{k}|_{\Omega_p} \quad \mbox{in } \Omega_p\,,\\
\lambda w_1^{k}-w_2^{k} &= f_2^{k} \quad \mbox{in } \Omega_p\,,\\
w_1^{k} = \frac{\partial w_1^{k}}{\partial \nu} & = 0  \quad \mbox{on } \partial\Omega_p \,. \end{aligned}
\end{eqnarray}
\item Update $w_2^{k} = w_2^{k-1}$.
\end{enumerate}
\textbf{end while}
\caption{Working algorithm}
\label{alg1}
\end{algorithm}

\begin{table}[!htbp] 
\centering 
\begin{tabular}{l l l l l} 
\toprule 
 $h$ & $\textbf{L}^2$(Velocity) & $\textbf{H}^1$(Velocity) & $L^2$(pressure)
\\ [0.5ex]
\midrule
1/4 & $1.17e-05$ &$3.29e-04$&$8.82e-05$\\ [0.5ex]
 1/6 & $5.05e-06[2.06]$ &$2.13e-04[1.07]$&$2.94e-05[2.71]$ \\ [0.5ex]
 1/8 & $2.53e-06[2.40]$ &$1.43e-04[1.38]$&$1.24e-05[2.99]$\\ [0.5ex]
 1/10 & $1.41e-06[2.63]$ &$1.00e-04[1.61]$&$6.20e-06[3.12]$\\ [0.5ex]
1/12 & $8.54e-07[2.75]$ &$7.29e-05[1.74]$&$3.44e-06[3.23]$\\ [0.5ex]
\bottomrule
\end{tabular}
\caption{Errors and convergence rates for $\bu$ and $p$ functions for different values of $h$.} 
\label{3dtab:2nd}
\end{table}

\begin{table}[!htbp] 
\centering 
\begin{tabular}{l l l l} 
\toprule 
 $h$ & $L^2$($w_1$-function) & $H^1$($w_1$-function) & $H^2$($w_1$-function)
\\ [0.5ex]
\midrule
1/4 & $8.54e-07$ &$6.79e-06$&$9.28e-05$\\ [0.5ex]
 1/6 & $5.16e-07[1.24]$ &$4.14e-06[1.22]$&$6.31e-05[0.95]$\\ [0.5ex]
 1/8 & $2.73e-07[2.21]$ &$2.33e-06[2.00]$&$4.76e-05[0.97]$\\ [0.5ex]
 1/10 & $1.82e-07[1.83]$ &$1.58e-06[1.74]$&$4.05e-05[0.74]$\\ [0.5ex]
1/12 & $1.16e-07[2.45]$ &$1.04e-06[2.32]$&$3.10e-05[1.46]$\\ [0.5ex]
\bottomrule
\end{tabular}
\caption{Errors and convergence rates for $w_1$ function for different values of $h$.} 
\label{3dtab:2ndw}
\end{table}

\clearpage
\newpage

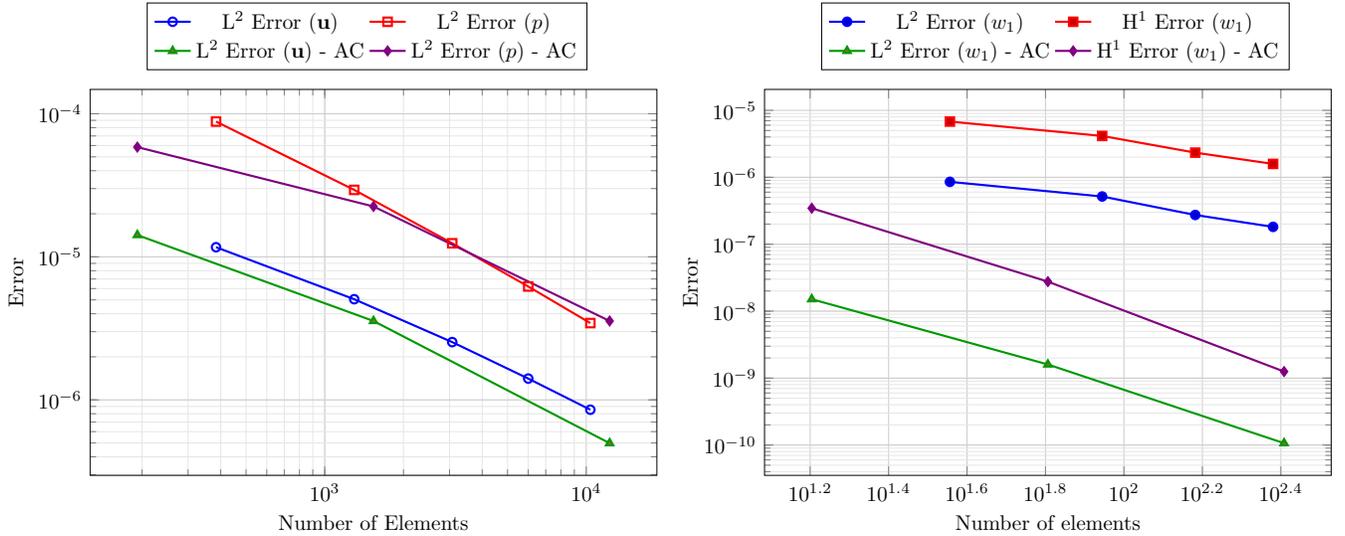
\begin{figure*}[h!]
    \centering
    \begin{subfigure}[t]{0.5\textwidth}
        \centering
\begin{tikzpicture}[scale=0.8]
\begin{loglogaxis}[
    width=11cm,
    height=8cm,
    xlabel={Number of Elements},
    ylabel={Error},
    legend style={at={(0.1,1.05)}, anchor=south west, legend columns=2},
    grid=both,
    major grid style={line width=.2pt,draw=gray!40},
    minor grid style={line width=.1pt,draw=gray!20},
    mark options={solid},
    every axis plot/.append style={line width=1pt},
    font=\small,
]

\addplot+[mark=o, blue] coordinates {
    (384, 1.16733e-05)
    (1296, 5.05498e-06)
    (3072, 2.53476e-06)
    (6000, 1.40912e-06)
    (10368, 8.53512e-07)
};
\addlegendentry{L$^2$ Error ($\bu$)}

\addplot+[mark=square, red] coordinates {
    (384, 8.82124e-05)
    (1296, 2.9376e-05)
    (3072, 1.2438e-05)
    (6000, 6.20411e-06)
    (10368, 3.4428e-06)
};
\addlegendentry{L$^2$ Error ($p$)}

\addplot+[mark=triangle*, green!60!black] coordinates {
    (192, 1.42e-05)
    (1536, 3.56e-06)
    (12288, 4.98e-07)
};
\addlegendentry{L$^2$ Error ($\bu$) - AC}

\addplot+[mark=diamond*, violet] coordinates {
    (192, 5.85e-05)
    (1536, 2.25e-05)
    (12288, 3.56e-06)
};
\addlegendentry{L$^2$ Error ($p$) - AC}

\end{loglogaxis}
\end{tikzpicture}
\end{subfigure}%
    ~ 
    \begin{subfigure}[t]{0.5\textwidth}
        \centering

\begin{tikzpicture}[scale=0.8]
\begin{loglogaxis}[
    width=11cm,
    height=8cm,
    xlabel={Number of elements},
    ylabel={Error},
    legend style={at={(0.1,1.05)}, anchor=south west, font=\small,legend columns=2},
    grid=both,
    major grid style={line width=.2pt,draw=gray!40},
    minor grid style={line width=.1pt,draw=gray!20},
    mark options={solid},
    every axis plot/.append style={line width=1pt},
    font=\small,
]

\addplot+[mark=*, blue] coordinates {
    (36, 8.54494e-07)
    (88, 5.16477e-07)
    (152, 2.7344e-07)
    (240, 1.81582e-07)
};
\addlegendentry{L$^2$ Error ($w_1$)}

\addplot+[mark=square*, red] coordinates {
    (36, 6.79267e-06)
    (88, 4.14476e-06)
    (152, 2.33341e-06)
    (240, 1.5842e-06)
};
\addlegendentry{H$^1$ Error ($w_1$)}

\addplot+[mark=triangle*, green!60!black] coordinates {
    (16, 1.509e-08)
    (64, 1.598e-09)
    (256, 1.066e-10)
};
\addlegendentry{L$^2$ Error ($w_1$) - AC}

\addplot+[mark=diamond*, violet] coordinates {
    (16, 3.45e-07)
    (64, 2.761e-08)
    (256, 1.253e-09)
};
\addlegendentry{H$^1$ Error ($w_1$) - AC}

\end{loglogaxis}
\end{tikzpicture}

\end{subfigure}
    \caption{A comparision of the errors against the number of elements with the results reported in \cite{A-Clark}.}
\label{fig:Comp}
\end{figure*}

In \Cref{fig:Comp}, we compare the errors against the number of elements with the results presented in \cite{A-Clark}. The results from \cite{A-Clark} are labeled as ``AC" in the legend. The errors for the fluid variables are comparable, as both approaches utilize Taylor–Hood elements for the fluid subproblem. However, as expected, the errors for the plate displacement are larger than the results reported in \cite{A-Clark}, as we use the lower order non-conforming P2-Morley elements for $w_1$ instead of P5 as in \cite{A-Clark}. Hence, the errors and the convergence rates are not directly comparable. The P2-Morley elements, being nonconforming, allow for simpler basis function construction and fewer degrees of freedom per element, making them more efficient and practical for large-scale or coupled fluid–structure interaction problems without severely compromising accuracy.

\clearpage
\newpage

\section{Conclusion}\label{sec:Conclude}
We developed a nonstandard mixed variational formulation for a fluid–structure interaction model coupling a 3D incompressible fluid with a 2D elastic plate, and established strongly continuous semigroup wellposedness through a novel approach. Our approach, based on introducing a Lagrange multiplier characterized by the boundary trace of the pressure, avoids the nonlocal maps of earlier work and provides simpler and efficient mixed finite element method for approximating solutions to the fluid-plate interaction systems. For the numerical approximation, we used Taylor-Hood elements for the fluid and $P2$–Morley elements for the plate problem within a domain decomposition scheme and presented numerical results to demonstrate the effectiveness of our approach. Further, we compared our errors to \cite{A-Clark}. Our method yields comparable fluid errors, while plate errors are larger due to the use of lower-order nonconforming elements. Nevertheless, the Morley elements significantly reduce computational cost and complexity, making them well-suited for large-scale coupled problems. This combination of a simpler theoretical framework and a practical numerical scheme advances the analysis and computation of fluid–plate interaction systems, with potential applications in aeroelasticity, fluid dynamics, and biomedicine.

\section{Acknowledgement}

\noindent 
The author, George Avalos, was partially supported by the NSF Grant DMS-1948942 and Simons Foundation Grant MP-TSM-00002793.\\
The author, Pelin G. Geredeli, was partially supported by the NSF Grant DMS-2348312. \\
The author, Hyesuk Lee,  was partially supported by the NSF Grants DMS-2207971 and DMS-2513073.

\end{document}